\documentclass[a4paper,12pt,openright]{report}
\usepackage{t1enc}
\usepackage[brazil]{babel}
\usepackage[latin1]{inputenc}
\usepackage{latexsym}
\usepackage{amsfonts}
\usepackage{amssymb}
\usepackage{amsmath}
\usepackage{multicol}
\usepackage[all]{xy}
\usepackage{indentfirst}

\newcommand{\sumi}{\overset{N}{\underset{i=1}{\sum}}}
\newcommand{\sumk}{\overset{n}{\underset{k=1}{\sum}}}
\newcommand{\sumj}{\overset{M}{\underset{j=1}{\sum}}}

\newcommand{\Dem}{{\bf Demonstração: }}
\newcommand{\fimdem}{\begin{flushright} $\blacksquare$ \end{flushright}}
\newtheorem{Def}{Definição}[chapter]
\newtheorem{Prop}[Def]{Proposição}
\newtheorem{Cor}[Def]{Corolário}
\newtheorem{Lem}[Def]{Lema}
\newtheorem{Teo}[Def]{Teorema}
\newtheorem{Ex}[Def]{Exemplo}
\newtheorem{Obs}[Def]{Observação}

\addto\captionsbrazil{
}

\begin{document}
\thispagestyle{empty}

\vspace{6.0 cm}

\centerline{ \large{\bf O caráter de Chern-Connes para}}
\centerline{ \large{\bf C$^*$-sistemas dinâmicos calculado}}
\centerline{ \large{\bf em algumas álgebras de }}
\centerline{ \large{\bf operadores pseudodiferenciais}}

\vspace{1.5 cm}

\centerline{ \large{David Pires Dias}}

\vspace{3.0 cm}

\centerline{ {\sc Tese Apresentada}}
\centerline{ {\sc ao}}
\centerline{ {\sc Instituto de Matemática e Estatística}}
\centerline{ {\sc da}}
\centerline{ {\sc Universidade de São Paulo}}
\centerline{ {\sc para}}
\centerline{ {\sc obtenção do grau}}
\centerline{ {\sc de}}
\centerline{ {\sc Doutor em Ciências}}

\vspace{2.5 cm}

\centerline{ {Área de Concentração: {\bf Matemática}}}
\centerline{ {Orientador: {\bf Prof. Dr. Severino Toscano do Rego Melo}}}

\vspace{2.0 cm}

\centerline{ - São Paulo, abril de 2008. - }

\newpage
\thispagestyle{empty}
\ \ \ 

\newpage
\thispagestyle{empty}

\vspace{6.0 cm}

\centerline{ \large{\bf O caráter de Chern-Connes para}}
\centerline{ \large{\bf C$^*$-sistemas dinâmicos calculado}}
\centerline{ \large{\bf em algumas álgebras de }}
\centerline{ \large{\bf operadores pseudodiferenciais}}

\vspace{3.0 cm}

\begin{flushright}
{\sf Este exemplar corresponde à redação }\\
{\sf final da tese \ devidamente \ corrigida }\\
{\sf e defendida por David \ Pires \ Dias e }\\
{\sf aprovada \ pela \ comissão \ julgadora.}
\end{flushright}

\vspace{3.0 cm}

\begin{flushright} {\it São Paulo, 11 de abril de 2008.} \end{flushright}

\vspace{2.0 cm}

Banca Examinadora:
\begin{itemize}
\item Prof. Dr. Severino Toscano do Rego Melo (Presidente) \ \  - \ \  IME-USP
\item Prof. Dr. Ricardo Bianconi \ \ - \ \ IME-USP
\item Prof. Dr. Carlos Eduardo Duran Fernandez \ \ - \ \ UNICAMP
\item Prof. Dr. Antônio Roberto da Silva \ \ - \ \ UFRJ
\item Prof. Dr. Fernando Raul Abadie Vicens \ \ - \ \ UR
\end{itemize}

\newpage
\thispagestyle{empty}
\ \

\newpage
\thispagestyle{empty}
.
\vspace{12cm}

\begin{flushright}
{\it Aos meus pais,  }\\
{\it Carlos e Maria Helena. }
\end{flushright}

\newpage
\thispagestyle{empty}
\ \

\newpage

\pagenumbering{roman}

\chapter*{Resumo}

Dado um C$^*$-sistema dinâmico $(A, G, \alpha)$ define-se um homomorfismo, denominado de caráter de Chern-Connes, que leva elementos de $K_0(A) \oplus K_1(A)$, grupos de K-teoria da C$^*$-álgebra $A$, em $H_{\mathbb{R}}^*(G)$, anel da cohomologia real de deRham do grupo de Lie $G$. Utilizando essa definição, nós calculamos explicitamente esse homomorfismo para os exemplos $(\overline{\Psi_{cl}^0(S^1)}, S^1, \alpha)$ e $(\overline{\Psi_{cl}^0(S^2)}, SO(3), \alpha)$, onde $\overline{\Psi_{cl}^0(M)}$ denota a C$^*$-álgebra gerada pelos operadores pseudodiferenciais clássicos de ordem zero da variedade $M$ e $\alpha$ a ação de conjugação pela representação regular (translações).

\ \ \ 
\newpage
\thispagestyle{empty}
\ \ \ 
\newpage

\chapter*{Abstract}

Given a C$^*$-dynamical system $(A, G, \alpha)$ one defines a homomorphism, called the Chern-Connes character,  that take an element in $K_0(A) \oplus K_1(A)$, the K-theory groups of the C$^*$-algebra $A$, and maps it into $H_{\mathbb{R}}^*(G)$, the real deRham cohomology ring of $G$. We explictly compute this homomorphism for the examples $(\overline{\Psi_{cl}^0(S^1)}, S^1, \alpha)$ and $(\overline{\Psi_{cl}^0(S^2)}, SO(3), \alpha)$, where $\overline{\Psi_{cl}^0(M)}$ denotes the C$^*$-álgebra gene\-rated by the classical pseudodifferential operators of zero order in the manifold $M$ and $\alpha$ the action of conjugation by the regular representation (translations).

\ \ \ 
\newpage
\thispagestyle{empty}
\ \ \ 
\newpage

\tableofcontents

\ \ \ 

\newpage

\pagenumbering{arabic}
\chapter*{Introdução}
\addcontentsline{toc}{chapter}{Introdução}
O propósito deste trabalho foi, desde o princípio, o de compreender o caráter de Chern-Connes, apresentado pela primeira vez por Alain Connes em \cite{AC1}, e posteriormente calculá-lo na C$^*$-álgebra gerada pelos operadores pseudodiferenciais clássicos de ordem zero da esfera, que designamos por $\overline{\Psi_{cl}^0(S^2)}$.

Por esse motivo o primeiro capítulo trata da construção (definição) do homomorfismo de Chern-Connes para um  C$^*$-sistema dinâmico $( A, G, \alpha)$, onde $A$ é uma C$^*$-álgebra, $G$ é um grupo de Lie e $\alpha$ é um homomorfismo contínuo de $G$ no grupo dos automorfismos de $A$, designado por $Aut(A)$, equipado com a topologia da convergência pontual. 

Grande parte das idéias desse capítulo baseiam-se em \cite{AC1} que é um artigo extremamente conciso. Por esse motivo se fez necessário detalhar muitas das passagens omitidas nesse artigo e para tal foi necessário recorrer a idéias de outros textos que não aparecem exatamente no mesmo contexto em que estamos trabalhando. Vale destacar alguns desses textos como  \cite{AC2} do mesmo autor, \cite{LO} mais voltado à parte algébrica, \cite{MK} e \cite{Fe} dentre outros. 
Além desses textos, alguns dos resultados que provamos neste primeiro capítulo foram demonstrados graças a conversas com colegas como Johannes Aastrup, Ricardo Bianconi e Bertrand Monthubert.


O segundo capítulo nasceu da necessidade de se trabalhar um pouco com a definição dada no primeiro capítulo para casos mais simples do que o já mencionado no primeiro parágrafo. Por este motivo o capítulo em questão apresenta o cálculo explícito do caráter de Chern-Connes para o C$^*$-sistema dinâmico em que a C$^*$-álgebra é $\overline{\Psi_{cl}^0(S^1)}$, isto é, a C$^*$-álgebra gerada pelos operadores pseudodiferenciais clássicos de ordem zero da variedade $S^1$. 

A beleza do segundo capítulo está no fato de se ver inteiramente a construção do cárater de Chern-Connes, feita no capítulo anterior, sendo aplicada num caso concreto. Destacam-se aqui vários fatos como o de $\Psi_{cl}^0(S^1)$ ser invariante pelo cálculo funcional holomorfo o que faz com que este exemplo se encaixe perfeitamente à definição utilizada; a possibilidade de mostrar de forma simples que a aplicação do índice $\delta_1$ é sobrejetora, pois neste caso apresentamos de forma clara o operador que é levado no gerador do contra-domínio; e também o fato do traço neste caso ser uma combinação linear \emph{real} de duas integrais para que as hipóteses sejam todas satisfeitas. Mas talvez o fato que mereça maior atenção seja o de que pudemos provar, utilizando os isomorfismos de K-teoria e também o auxílio computacional\footnote{Para os cálculos em questão foi utilizado o "software" \ Maple 6.}, o \emph{caso particular} (\emph{a}) enunciado por A. Connes em \cite{AC1}. 

O terceiro capítulo apresenta o cálculo do caráter de Chern-Connes para a C$^*$-álgebra gerada pelos operadores pseudodiferenciais clássicos de ordem zero da esfera, denotada por $\overline{\Psi_{cl}^0(S^2)}$, que foi nosso objetivo desde o princípio do trabalho. Apresentamos na primeira parte desse capítulo os grupos de K-teoria de $C(SS^2)$, onde $SS^2 \subset TS^2$ denota o sub-fibrado das esferas unitárias do fibrado tangente de $S^2$ que denominamos fibrado das esferas de $S^2$, e estes foram obtidos utilizando-se Mayer-Vietoris para K-teoria de C$^*$-álgebras, ferramenta esta apresentada como em \cite{SM}, mas que também aparece de forma mais geral em \cite{BB}. Cabe ressaltar que os resultados aqui obtidos, isto é, que $K_1(C(SS^2)) = \mathbb{Z}$ e que $K_0(C(SS^2))= \mathbb{Z} \oplus \mathbb{Z}_2$, coincidiram com os resultados ainda não publicados por F. Rochon em \cite{FR}, mas vale dizer que nossos resultados foram encontrados através de ferramentas diferentes das utilizadas em \cite{FR}.

Já na segunda seção deste terceiro capítulo, utilizando os resultados obtidos na seção anterior, calculamos os grupos de K-teoria de $\overline{\Psi_{cl}^0(S^2)}$. E de posse dos grupos de K-teoria da C$^*$-álgebra $\overline{\Psi_{cl}^0(S^2)}$ pudemos chegar ao nosso objetivo principal e apresentar o caráter de Chern-Connes, como desejado desde o início.

O quarto capítulo apresenta o cálculo da K-teoria do fibrado das coesferas do toro, pois planejo, em seguida, utilizar estes grupos para calcular o caráter de Chern-Connes em $\overline{\Psi_{cl}^0(T)}$, ou seja, a C$^*$-álgebra gerada pelo operadores pseudodiferenciais clássicos de ordem zero do Toro. Mas para resolver este problema, em que pretendo continuar trabalhando, talvez seja necessário pensar em alguns outros, como por exemplo, uma fórmula mais geral para o caráter de Chern-Connes quando calculado em $K_1$, assim como fizeram Charlotte Wahl em \cite{CW} e também Ezra Getzler em \cite{EG} para outras generalizações do caráter de Chern-Connes.  Observe que tal fórmula não foi necessária nos dois exemplos que demos anteriormente, já que no primeiro, isto é, em $\overline{\Psi_{cl}^0(S^1)}$, devido a dimensão do grupo de Lie $S^1$, precisamos apenas da fórmula para o primeiro termo ímpar do somatório e no segundo exemplo (em $\overline{\Psi_{cl}^0(S^2)}$) não se fez necessária, pois $K_1(\overline{\Psi_{cl}^0(S^2)}) = 0$.

Vale lembrar que em todos os exemplos citados a C$^*$-álgebra em questão, isto é, a C$^*$-álgebra gerada pelos operadores pseudodiferenciais clássicos de ordem zero das variedades $S^1$, $S^2$ e $T$ contém a álgebra dos operadores compactos $\mathcal{K}$ e seus comutadores são compactos. Utilizando-se este fato e os resultados de "comparison algebras" \ de H. Cordes \cite{CO}, que na verdade já eram conhecidos por Konh e Nirenberg \cite{KN} e também por Gohberg \cite{GO} e Seeley \cite{Se}, temos o quociente da álgebra por $\mathcal{K}$ comutativo o que faz com que nossos exemplos não fiquem muito longe do caso clássico comutativo. Pretendemos trabalhar no cálculo explícito do homomorfismo de Chern-Connes para algumas álgebras geradas por operadores pseudodiferenciais que não possuem esta propriedade, vide por exemplo a álgebra estudada em \cite{TO1}. 

Além desses capítulos esta tese apresenta logo a seguir uma lista de convenções e notações iniciais. E na tentativa de auxiliar àqueles menos familiarizados com K-teoria e/ou com alguns resultados sobre grupos de Lie, esta tese também dispõe de dois apêndices. O primeiro versando sobre alguns resultados, necessários para a compreensão do cárater de Chern-Connes no contexto de C$^*$-sistemas dinâmicos, sobre representações em grupos de Lie. Cabe aqui agradecer a Daniel V. Tausk por sua inestimável cooperação e esclarecimentos sobre diversos tópicos desse apêndice. Existe ainda um segundo apêndice sobre K-teoria, ferramenta essencial para o entendimento e os cálculos do homomorfismo de Chern-Connes tema central deste trabalho.

\ \ \

\newpage

\ \ \

\newpage

\addcontentsline{toc}{chapter}{Notações e convenções}
\chapter*{Notações e convenções}
\begin{itemize}

\item $S^1 = \{ z \in \mathbb{C} : |z| = 1\}$ \ \ \ , \ \ \ $T = S^1 \times S^1$ \ \ \ e \ \ \ $D = \{ z \in \mathbb{C}: |z| \leqslant 1 \}$.

\item $M_a$ é o operador de multiplicação por $a \in C^\infty (S^1)$ e $\mathcal{L}(\mathcal{H})$ o conjunto dos operadores limitados de $\mathcal{H}$ como em \ref{BH}. 

\item $F_d$ designa a transformada de Fourier discreta, isto é, $F_d: L^2 (S^1) \rightarrow l^2(\mathbb{Z})$, onde $(F_du)_j = \frac{1}{\sqrt{2 \pi}} \int_{- \pi}^\pi \tilde{u}(\theta) e^{-ij \theta} d \theta$, $j \in \mathbb{Z}$ e $\tilde{u}(\theta) = u (e^{i \theta})$.

\item $S^*M$ denotará o fibrado das coesferas de uma variedade riemanniana $M$ de dimensão $n$, isto é, o conjunto de todos os pontos $(x, \xi ) \in T^*M$, fibrado cotangente de $M$,  tais que $\underset{i, j = 1}{\overset{n}{\sum}} m_{ij}(x) \xi_i \xi_j = 1$, onde $m$ é a métrica riemanniana.

\item $CS(\mathbb{Z}) = \{ (a_i)_{i \in \mathbb{Z}} \ : \ \underset{i \rightarrow \infty}{\lim} a_i = a(\infty) \ \ \textrm{e} \ \  \underset{i \rightarrow  -\infty}{\lim} a_i = a(-\infty) \ \  \textrm{existem} \}$.

\item $a(D_\theta) = F_d^{-1}M_aF_d$, onde $M_a$ é o operador de multiplicação pela seqüência $(a_i) \in CS(\mathbb{Z})$.

\item Os símbolos $\mathfrak{1}$, $\mathfrak{0}$ e $\mathfrak{z}$ serão utilizados para designar, respectivamente, as aplicações $\mathfrak{1} (z) = 1$, $\mathfrak{0}(z) = 0$ e $\mathfrak{z} (z) = z$, $z \in S^1$, como pode ser visto em \ref{cded} e \ref{Ks1}.

\item O conjunto $SA$ é a suspensão da C$^*$-álgebra $A$, isto é, 
\begin{equation*}
\begin{split}
SA &= \{ f \in C([0,1],A): f(0) = f(1)= 0 \} \\
 & = \{ f \in C([0,2 \pi],A): f(0) = f(2 \pi)= 0 \} \\ 
 & = \{ f \in C(S^1,A): f(1) = 0 \} \ \ = \ \ C_0(]0,1[, A)\\
\end{split}
\end{equation*}
como em \ref{SA}.

\item  Sendo $A$ uma $C^*$-álgebra as aplicações $\theta_A$ e $\beta_A$ são os isomorfismos da K-teoria complexa dados por $\theta_A : K_1(A) \rightarrow K_0(SA)$ e  $\beta_A : K_0(A) \rightarrow K_1(SA)$ respectivamente, para maiores detalhes vide \ref{thetaa} e \ref{Bott}.

\end{itemize}

\ \ \

\newpage

\chapter{O caráter de Chern-Connes para C$^*$-sistemas dinâmicos}
O objetivo deste capítulo é o de apresentar com mais detalhes a definição do caráter de Chern dada por Alain Connes em \cite{AC1}. Nesse artigo Connes generaliza a idéia de caráter de Chern, já consagrada em geometria diferencial.

Na geometria diferencial o caráter de Chern nada mais é do que um morfismo que relaciona elementos de um K-grupo a elementos de um grupo de cohomologia (de deRham). Acontece porém, que tal caráter nasceu naturalmente para cálculos em que os K-grupos provinham de Álgebras de Banach comutativas e Alain Connes estendeu estes resultados e cálculos, já clássicos para o caso comutativo, para C$^*$-sistemas dinâmicos (caso este em que existe uma ação de grupo envolvida e que estudaremos neste capítulo) e também para álgebras de Banach não-comutativas utilizando teorias mais sofisticadas como a da cohomologia cíclica, cíclica periódica, etc. (casos estes que podem ser encontrados em \cite{AC2} e \cite{LO}).

Vale lembrar que a definição mais usual do carácter de Chern da geometria diferencial é dada através do cálculo do traço da exponencial de uma curvatura do fibrado vetorial de uma variedade diferenciável e que tal traço é um elemento do grupo de cohomologia de deRham desta variedade. O que apresentaremos a seguir é um detalhamento da definição do caráter de Chern para C$^*$-sistemas dinâmicos e isto será feito de maneira análoga a definição mais usual deste homomorfismo.

Como já dito a apresentação a seguir é baseada na extensão do caráter de Chern introduzida por Alain Connes em \cite{AC1} e por este motivo passaremos a denominar tal homomorfismo por caráter de Chern-Connes.

\section{O Caráter de Chern-Connes}

Apresentaremos nesta seção a definição do caráter de Chern-Connes introduzido por Alain Connes em seu artigo \cite{AC1} preenchendo também alguns detalhes omitidos pelo autor no referido artigo. 

Seja $(A, G, \alpha)$ um C$^*$-sistema dinâmico, isto é, $A$ é uma C$^*$-álgebra unital, $G$ um grupo de Lie e $\alpha:G \rightarrow Aut(A)$ uma ação contínua na topologia forte de operadores, isto é, para todo $a \in A$ a aplicação $\alpha(a): g \mapsto \alpha_g(a)$ é contínua na norma.

Diremos que $a \in A$ é de classe $C^k$ ($1 \leqslant k \leqslant \infty$)  se, e somente se, a aplicação $\alpha(a)$ é de classe $C^k$. A álgebra involutiva $A^\infty = \{ a \in A : a  \ \textrm{é de classe} \ C^\infty \}$ é densa (na norma) em $A$. Este resultado é conhecido como Teorema de G\aa rding e sua demonstração pode ser encontrada em \ref{Garding}.

Cabe observar que $A^\infty$ possui uma estrutura de *-álgebra de Fréchet induzida pela aplicação $\alpha :A^\infty \rightarrow C^\infty(G,A)$, dada por $a \mapsto \alpha(a)$. Esta topologia, induzida pelas seminormas canônicas de $C^\infty (G,A)$, torna a inclusão $A^\infty \hookrightarrow A$ contínua. De fato, como $\alpha$ é um automorfismo, temos $$|| \alpha(a_i)||_\infty = \underset{g \in G}{\sup}|| \alpha_g (a_i)|| = ||a_i ||.$$

Além disso, se $a \in A^\infty \subset A$ é invertível em $A$, então $a^{-1} \in A^\infty$, isto porque $\alpha_g(a^{-1}) = \alpha_g(a)^{-1}$. Assim o conjunto do invertíveis de $A^\infty$ é aberto e portanto a inversão é contínua, vide o corolário da página 115 de \cite{LW}. 

Como a operação de inversão é contínua, então as integrais de Cauchy que nos dão o cálculo funcional holomorfo convergem em $A^\infty$, em outras palavras, $A^\infty$ é invariante pelo cálculo funcional holomorfo e portanto a inclusão de $A^\infty$ em $A$ induz os isomorfismos $K_0(A^\infty) \simeq K_0(A)$ e $K_1(A^\infty) \simeq K_1(A)$. Este é um resultado bastante conhecido de K-teoria, que afirma que dadas uma C$^*$-álgebra $A$ e uma subalgebra $A^\infty$ densa em $A$ e invariante pelo cálculo funcional holomorfo, então $K_i(A^\infty) \simeq K_i(A)$, para $i=0$ ou $1$. A demonstração deste fato pode ser vista com mais detalhes no apêndice 3 de \cite{AC3} ou também como enunciado nas seções $5$ e $8$ de \cite{BB}.

Nossa intenção nesta seção é a de definir o caráter de Chern-Connes que, como veremos a seguir, é um homomorfismo entre os grupos de K-teoria de uma C$^*$-álgebra e os de cohomologia de um grupo de Lie. Podemos então, utilizando os isomorfismos vistos no parágrafo anterior, isto é, $K_0(A^\infty) \simeq K_0(A)$ e $K_1(A^\infty) \simeq K_1(A)$, fazer a opção de utilizar um m.p.f.g.\footnote{módulo à direita, exceto quando mencionado o contrário, projetivo e finitamente gerado}  $M^\infty$ sobre $A^\infty$, ao invés de $M$ um m.p.f.g. sobre $A$. Tal opção não trará diferença à definição do caráter de Chern-Connes, no entanto, nos permitirá derivar os elementos da álgebra, já que estes estarão sempre em $A^\infty$. 



Como todo m.p.f.g é um somando direto de um módulo livre com base finita, vide proposição $3.10$ de \cite{JA} e seu corolário, podemos então afirmar que existe um idempotente\footnote{mais precisamente uma família de idempotentes.} $e \in M_n(A^\infty)$, para algum $n \in \mathbb{N}$, tal que $e((A^\infty)^n)$ e $M^\infty$ são isomorfos. Por $4.6.2$ de \cite{BB} sabemos que dentre os idempotentes desta mesma classe podemos tomar um que satisfaça as condições
 $$e^2 = e \ \ \ \textrm{ e } \ \ \ e^*=e.$$ 
Ou seja, $e$ é uma projeção que passaremos a designar por $p \in M_n(A)$. 

\begin{Def} \label{Der}
Definimos $\delta$ uma representação da álgebra de Lie do grupo $G$, de\-signada por $\mathfrak{g}$, na álgebra de Lie das derivações de $A^\infty$, dada por 
$$ \delta_X (a) = \underset{t \rightarrow 0}{\lim} \frac{\alpha_{g_t} (a) - a}{t} \ , \ \ \textrm{ onde } \ \  g_0 ' = X \ \ \textrm{ e } \ \ a \in A^\infty. $$
\end{Def}

A demonstração de que $\delta$ é de fato uma representação, isto é, de que a definição anterior faz sentido, pode ser vista em \ref{representacao}. 

Algumas vezes, ao invés de nos depararmos com $a \in A$  {\bf(} $a \in A^\infty$ {\bf)}, teremos $a\in M_n(A)$ {\bf(} $a \in M_n(A^\infty)$ {\bf)} e a extensão da definição das aplicações $\alpha_g$ {\bf(} $\delta_X$ {\bf)} de $A$ {\bf(} $A^\infty$ {\bf)}, para $M_n(A)$ {\bf(} $M_n(A^\infty)$ {\bf)} é feita pontualmente, ou melhor, em cada entrada da matriz, isto é, para $a=(a_{ij}) \in M_n(A)$, temos $\alpha_g(a) = (\alpha_g(a_{ij}))$ {\bf(} para $a=(a_{ij}) \in M_n(A^\infty)$, temos $\delta_X(a) = (\delta_X(a_{ij}))$ {\bf)}.

Com o que foi apresentado até o momento, nos encontramos aptos a definir o que é uma conexão $\nabla$ sobre o m.p.f.g. $M^\infty$. 

\begin{Def}
Dado o m.p.f.g.  $M^\infty$ sobre $A^\infty$ uma conexão sobre $M^\infty$ é uma aplicação linear $\nabla: M^\infty \rightarrow M^\infty \otimes \mathfrak{g}^*$, que satisfaz  
$$ \nabla_X(\xi a) = \nabla_X(\xi)a + \xi \delta_X (a)  \ , \ \ \forall \xi \in M^\infty , \ \ \forall X \in \mathfrak{g} \ \ \ \textrm{e} \ \ \ \forall a \in A^\infty.$$
\end{Def}


Para todo $M^\infty \simeq p ((A^\infty)^n)$ m.p.f.g. sobre $A^\infty$, com $p $ como acima, existe uma conexão 
chamada de conexão grassmanianna ou de Levi-Civita (em analogia à definição usual da conexão de Levi-Civita da geometria diferencial). Tal conexão é dada por:
$$ \nabla_X^0(\xi) = p \delta_X(\xi) \ \ \in \ \ M^\infty, \ \ \ \ \forall \xi \in M^\infty \ \textrm{e} \ X \in \mathfrak{g}.$$

De fato, a linearidade de $\nabla_X^0$ decorre da linearidade de $\delta_X$. Através do isomorfismo $p \big( (A^\infty)^n \big) \simeq M^\infty$ e da igualdade $p^2 = p$, dados $X \in \mathfrak{g}$, $a \in A^\infty$ e $\xi \in M^\infty$, temos 
$$\nabla_X^0 (\xi a) = p \delta_X(\xi a) = p \delta_X(\xi)a + p \xi \delta_X (a) = \nabla_X^0(\xi)a + \xi \delta_X(a),$$
 ou seja, $\nabla^0$ é conexão.

Através da representação $\delta$ temos o complexo $\Omega = A^\infty \otimes \Lambda \mathfrak{g}^*$ das formas dife\-renciáveis invariantes à esquerda sobre o grupo $G$ com coeficientes em $A^\infty$. Assim munimos $\Omega$ com uma estrutura de álgebra \footnote{Observe que esta estrutura não é necessariamente comutativa, isto porque $$w_1 \wedge w_2 \not= (-1)^{gr(w_1)gr(w_2)} w_2 \wedge w_1.$$}, dada pelo produto tensorial de $A^\infty$ com a álgebra exterior $\Lambda \mathfrak{g}^*$, cuja derivação exterior $d$ satisfaz: 
\begin{itemize}
\item $\delta_X(a) = da(X)$, $\forall a \in A^\infty$ e $X \in \mathfrak{g}$; 
\item $d(w_1 \wedge w_2) = d(w_1) \wedge w_2 + (-1)^k w_1 \wedge d(w_2)$, $\forall w_1 \in \Omega^k$ e $w_2 \in \Omega$ e 
\item $d^2 (w) =  0$, $\forall w \in \Omega$. 
\end{itemize}

Note que $A^\infty \simeq A^\infty \otimes {1} \subset \Omega$ e portanto $\Omega$ pode ser visto como um bimódulo sobre $A^\infty$. 

Além disso, podemos identificar $End(M^\infty) \simeq End(p (A^\infty)^n)$ com $p M_n(A^\infty) p$. Para isto basta tomar o isomorfismo 
\begin{equation*}
\begin{split}
 \phi : End(M^\infty) & \rightarrow p M_n(A^\infty) p \\
                       f & \mapsto pf(p)p .
\end{split}
\end{equation*}

\begin{Lem} \label{conexao}
Toda conexão $\nabla$ sobre $M^\infty$ é da forma 
$$\nabla_X(\xi) = \nabla_X^0(\xi) + \Gamma_X (\xi) \ , \ \ \forall \xi \in p \big( (A^\infty)^n \big) \ \ \textrm{ e } \ \ X \in \mathfrak{g},$$ com $\Gamma \in p M_n(\Omega^1) p$ unicamente determinada por $\nabla$.
\end{Lem}
\Dem 
Observe que $\gamma = \nabla - \nabla^0$ é $A^\infty$-linear , isto é, $\gamma_X (\xi a ) = \gamma_X(\xi) a$, para todo $X \in \mathfrak{g}$, $\xi \in M^\infty$ e $a \in A^\infty$. 

Utilizando o isomorfismo $\phi: End(M^\infty) \rightarrow p M_n(A^\infty) p$ podemos notar que existe $\Gamma_X \in M_n(A^\infty)$ tal que $\gamma_X(\xi) = p \Gamma_X p \xi$, $\forall \xi \in M^\infty$ e além disso, a aplicação $X \mapsto \Gamma_X$ é linear, isto é, $\Gamma \in p M_n(\Omega^1) p$.
\fimdem

\begin{Def}
Seja $\nabla$ uma conexão sobre o m.p.f.g. $M^\infty$, definimos a curvatura associada a $\nabla$ como sendo o elemento $\Theta \in End_{A^\infty}(M^\infty) \otimes \Lambda^2\mathfrak{g}^*$ dado por 
$$\Theta (X,Y) = \nabla_X\nabla_Y-\nabla_Y\nabla_X-\nabla_{[X,Y]}  \ , \ \ \ \forall X, \ Y \in \mathfrak{g}.$$
\end{Def}

\begin{Obs} \label{pdp}
Note que as curvaturas assim definidas são 2-formas e utilizando a identificação $End(M^\infty)$ com $pM_n(A^\infty) p \subset M_n(A^\infty)$ podemos mostrar que a curvatura asso\-ciada a conexão grassmanianna $\nabla^0$ é a 2-forma $ \Theta_0 = pdp \wedge dp \in \Omega^2$. 
\end{Obs}
Vejamos
\begin{equation*}
\begin{split}
\Theta_0(X,Y) &= p(\Theta_0(X,Y)(p))p = p(\nabla^0_X\nabla^0_Y(p)-\nabla^0_Y\nabla^0_X(p)-\nabla^0_{[X,Y]}(p)  )p \\
&= p(p\delta_X(p\delta_Y(p)) - p\delta_Y(p\delta_X(p))-p\delta_{[X,Y]}(p)  )p \\
\end{split}
\end{equation*}
e utilizando-se o fato de que $\delta$ é uma representação de álgebras de Lie, como provado em \ref{representacao}, temos
\begin{equation*} 
\begin{split}
\Theta_0(X,Y) &= p(\delta_X(p)\delta_Y(p)+\delta_X\delta_Y(p)-\delta_Y(p)\delta_X(p)-\delta_Y\delta_X(p) -\delta_X\delta_Y(p) + \delta_Y\delta_X(p) ) \\
&= p(\delta_X(p)\delta_Y(p) - \delta_Y(p)\delta_X(p)) = p dp \wedge dp (X,Y).
\end{split}
\end{equation*}

\begin{Prop} \label{curv}
Com a notação do Lema \ref{conexao} podemos mostrar que a curvatura associada a conexão $\nabla = \nabla^0 + \Gamma$ é dada por 
$$ \Theta = \Theta_0 + p (d \Gamma + \Gamma \wedge \Gamma)p .$$
\end{Prop}
\Dem
Seja $\nabla = \nabla^0 + \Gamma$ uma conexão como no Lema \ref{conexao}. Utilizando-se formas tensorizadas podemos provar, de forma análoga a feita nos itens \emph{(a)} e \emph{(f)} da proposição $2.25$ de \cite{Wa}, que   
$$d \Gamma (X, Y) = \delta_X \Gamma_Y - \delta_Y \Gamma_X - \Gamma_{[X,Y]}, $$ 
para quaisquer $X$ e $Y$ em $\mathfrak{g}$.

Além disso, $(\Gamma \wedge \Gamma)(X,Y) = \Gamma_X \Gamma_Y - \Gamma_Y \Gamma_X$.

Assim ao trocarmos $\nabla$ por $\nabla^0 + \Gamma$ na definição de curvatura temos para quaisquer $ X, \ Y \in \mathfrak{g}$ 
\begin{equation*}
\begin{split}
\Theta (X,Y) &= \nabla_X\nabla_Y-\nabla_Y\nabla_X-\nabla_{[X,Y]} \\
            &=  (\nabla^0_X + \Gamma_X)(\nabla^0_Y + \Gamma_Y) - (\nabla^0_Y + \Gamma_Y)(\nabla^0_X + \Gamma_X) - (\nabla^0_{[X,Y]} + \Gamma_{[X,Y]}) \\
           &= \nabla^0_X \nabla^0_Y + \nabla^0_X \Gamma_Y + \Gamma_X \nabla_Y^0 + \Gamma_X \Gamma_Y - \nabla^0_Y \nabla^0_X - \nabla^0_Y \Gamma_X + \\
      & \hspace{7.0 cm}    - \Gamma_Y \nabla_X^0 - \Gamma_Y \Gamma_X -  \nabla^0_{[X,Y]} - \Gamma_{[X,Y]} \\
\end{split}
\end{equation*}

Utilizando as definições da conexão de Levi-Civita $\nabla^0_X  = p \delta_X$ e também da curvatura associada a esta, isto é, $\Theta_0(X,Y) = \nabla^0_X\nabla^0_Y-\nabla^0_Y\nabla^0_X-\nabla^0_{[X,Y]}$ , temos

\begin{equation*}
\begin{split}
\Theta (X,Y) &= \Theta_0 (X,Y) + p \delta_X (\Gamma_Y) - p \delta_Y ( \Gamma_X ) + \Gamma_X p \delta_Y - \Gamma_Y p \delta_X + \\
& \hspace{8.5 cm} + \Gamma_X \Gamma_Y - \Gamma_Y \Gamma_X - \Gamma_{[X,Y]} \\
          &= \Theta_0 (X,Y) + p \delta_X \Gamma_Y + p \Gamma_Y \delta_X - p \delta_Y \Gamma_X - p \Gamma_X \delta_Y + \Gamma_X p \delta_Y + \Gamma_X p \delta_Y +\\ 
         &\hspace{8.5 cm}                  - \Gamma_{[X,Y]} +  \Gamma_X \Gamma_Y - \Gamma_Y \Gamma_X  \\
\end{split}
\end{equation*}

Lembrando que $\Gamma \in p M_n(\Omega^1)p$ e portanto $p \Gamma_X \delta_Y = \Gamma_X p \delta_Y$ e $ p \Gamma_Y \delta_X = \Gamma_Y p \delta_X$, temos
\begin{equation*}
\begin{split}
\Theta (X,Y) &= \Theta_0 (X,Y) + p ( \delta_X \Gamma_Y - \delta_X \Gamma_Y - \Gamma_{[X,Y]})  + p (\Gamma_X \Gamma_Y - \Gamma_Y \Gamma_X) \\ 
          &= \Theta_0 (X,Y) + p(d \Gamma)p (X,Y) + p(\Gamma \wedge \Gamma)p(X,Y) .
\end{split}
\end{equation*}

Ou seja, podemos descrever qualquer curvatura $\Theta$ por
$$ \Theta = \Theta_0 + p (d \Gamma + \Gamma \wedge \Gamma)p.$$
\fimdem




Utilizando agora o fato de $p \in M_n(A^\infty)$ ser projeção, ou mais precisamente $p^2= p$, podemos provar algumas identidades que nos serão úteis no decorrer desta seção.

\begin{Lem} \label{identidades}
Seja $p\in M_n(A^\infty)$ uma projeção e $d$ a derivação exterior, então:
\begin{enumerate}
\item $p dp p =0$;
\item $p dp \wedge dp = p (dp \wedge dp) p = (dp \wedge dp) p$;
\item $(2p-1)dp = -dp(2p-1)$.
\item $(p dp \wedge dp)^k = p (dp \wedge dp)^k$, $\forall k \in \mathbb{N}$;
\item $d(p dp \wedge dp) = dp \wedge (p dp \wedge dp) + (p dp \wedge dp) \wedge dp$.
\end{enumerate}
\end{Lem}
\Dem $(1)$ Como $p^2 = p$, então $d(p^2) = dp$ e pela regra de Leibniz 
$dp p + p dp = dp,$ multiplicando-se por $p$ ambos os lados da igualdade e utilizando novamente que $p^2 = p$, temos $pdpp = 0$.

$(2)$ Segue disto que $d(pdpp)= 0$, logo $dp \wedge dp p + p d^2p p - p dp \wedge dp = 0$ e como $d^2 = 0$, temos $dp \wedge dp p = p dp \wedge dp$. Utilizando o raciocínio anterior temos $p dp \wedge dp = p (dp \wedge dp) p = (dp \wedge dp) p$.

$(3)$ Pela demonstração de (1) podemos observar que 
\begin{equation*}
\begin{split}
(2p-1)dp &= 2pdp -dp = pdp + (pdp-dp) = p dp - dp p  \\
 &= (dp - dpp) - dpp = -2dpp+dp = -dp(2p-1).
\end{split}
\end{equation*}

$(4)$ Por (2) vemos que $p$ comuta com $dp \wedge dp$ logo $$(p dp \wedge dp)^k  = p^k (dp \wedge dp)^k = p (dp \wedge dp)^k.$$

$(5)$ Novamente por (2) temos $p dp \wedge dp = p dp \wedge dp p$, logo 
\begin{equation*}
\begin{split}
d(p dp \wedge dp) &= d (p dp \wedge dp p) = dp \wedge dp \wedge dp p + p dp \wedge dp \wedge dp \\
&= dp \wedge (p dp \wedge dp) + (p dp \wedge dp) \wedge dp.
\end{split}
\end{equation*}

\fimdem

A segunda afirmação deste lema também nos diz que $\Theta_0 = p\Theta_0 = \Theta_0 p = p\Theta_0p$ e a quinta que $d(\Theta_0) = dp \wedge \Theta_0 + \Theta_0 \wedge dp$. Assim 
\begin{equation*}
\begin{split}
d(\Theta_0) &= d (p \Theta_0 p) = dp \wedge \Theta_0 p + p d(\Theta_0)p + p \Theta_0 \wedge dp \\
            &= dp \wedge \Theta_0 + \Theta_0 \wedge dp + p d(\Theta_0)p \\
            &= d(\Theta_0) + p d(\Theta_0) p,
\end{split}
\end{equation*}
ou seja, $pd(\Theta_0)p = 0$.

\begin{Def} \label{tracog}
Uma aplicação $\tau: A \rightarrow \mathbb{C}$ é dita um traço finito $G$-invariante se para quaisquer elementos $a$ e $b$ de $A$ e $g \in G$, temos
\begin{enumerate}
\item $\tau$ é um funcional linear contínuo que satisfaz $\tau(ab) = \tau(ba)$ \ \ \ \ \ \ \ \ \ \ (traço);
\item $\tau$ é positivo e portanto $\tau(a^*) = \overline{\tau(a)}$ (vide os comentários após \ref{KdeC}) \ \ \ \ \ e 
\item $\tau(\alpha_g(a)) = \tau(a)$. 
\end{enumerate}
\end{Def}

Seja $\tau$ um traço finito $G$-invariante sobre $A$. Para todo $k \in \mathbb{N}$ existe uma única aplicação $k-$linear $\tau_k : \underbrace{\Omega \times \ldots \times \Omega}_{k} \rightarrow \Lambda \mathfrak{g}^*$, tal que $$\tau_k (a_1 \otimes w_1, \ldots, a_k \otimes w_k) = \tau (a_1 \ldots a_k)w_1 \wedge \ldots \wedge w_k.$$ 

Vale ressaltar que $\tau_0 = \tau$ é simplesmente o traço da álgebra $A$ e também que algumas vezes $\tau_k$ será entendido como a composição da aplicação $k-$linear definida acima com o traço usual das matrizes, já que estaremos trabalhando com elementos de $M_n(\Omega)$. 

Cabe também observar que $\tau_k (a_1 \otimes w_1, \ldots, a_k \otimes w_k) = \tau_1 (a_1 \ldots a_k \otimes w_1 \wedge \ldots \wedge w_k)$ e este fato será constantemente utilizado durante esta seção.

\begin{Lem} \label{traco}
A aplicação $\tau_k: \underbrace{\Omega \times \ldots \times \Omega}_{k} \rightarrow \Lambda \mathfrak{g}^*$ acima definida é um traço graduado, isto é, para $ i = 1, 2, \ldots , k-1 $ vale 
$$\tau_k(\theta_1, \ldots, \theta_k) = (-1)^{gr(\theta_i)gr(\theta_{i+1})} \tau_k(\theta_1, \ldots, \theta_{i-1},\theta_{i+1}, \theta_{i}, \theta_{i+2} \ldots \theta_k) .$$
\end{Lem}
\Dem
Provaremos primeiramente que $\tau_k$ está bem definido e isto se faz necessário, já que os elementos de $\Omega$ não possuem representação única.

Sejam $\sumi a_i \otimes w_i = \sumj b_j \otimes \lambda_j$ elementos de $\Omega$ e $\{ X_1, X_2, \ldots X_n\}$ uma base de $\mathfrak{g}$. Temos então
$$\sumk \left( \sumi w_i(X_k)a_i \right) \otimes \hat{X}_k  = \sumk  \left( \sumj \lambda_j(X_k)b_j \right) \otimes \hat{X}_k ,$$
logo 
$$\sumi w_i(X_k)a_i = \sumj \lambda_j(X_k)b_j.$$

Disto
\begin{equation*}
\begin{split}
\tau_1(\sumi a_i \otimes w_i) &= \sumi \tau(a_i)w_i =  \sumi \tau(a_i)\sumk w_i(X_k)\hat{X}_k \\
   &= \sumk \tau(\sumi w_i(X_k) a_i) \hat{X}_k = \sumk \tau(\sumj \lambda_j(X_k)b_j) \hat{X}_k \\
   &= \sumj \tau(b_j)\sumk \lambda_j(X_k)\hat{X}_k = \sumj \tau(b_j)\lambda_j \\
   &= \tau_1(\sumj b_j \otimes \lambda_j).
\end{split}
\end{equation*}

Provamos assim que $\tau_1$ está bem definido, mas decorre da observação anterior a este lema que isto basta para que $\tau_k$ esteja bem definida. 

Observemos agora que o sinal é alterado a cada permutação de elementos consecutivos em $\tau_k(a_1 \otimes w_1, \ldots, a_k \otimes w_k)$, de fato 
\begin{equation*}
\begin{split}
& \tau_k( a_1 \otimes w_1 , \ldots , a_i \otimes w_i, a_{i+1} \otimes w_{i+1}, \ldots a_k \otimes w_k) \\
&= \tau(a_1 \ldots a_i a_{i+1} \ldots a_k )w_1 \wedge \ldots \wedge w_{i-1} \wedge w_i \wedge w_{i+1} \ldots w_k \\
&= \tau(a_1 \ldots  a_{i-1}a_{i+1}a_ia_{i+2} \ldots a_k) w_1 \wedge \ldots \wedge w_{i-1} \wedge w_i \wedge w_{i+1} \wedge \ldots w_k \\
&= (-1)^{gr(w_i)gr(w_{i+1})} \tau(a_1 \ldots a_{i-1}a_{i+1}a_ia_{i+2} \ldots a_k) w_1  \ldots  w_{i-1}\wedge w_{i+1} \wedge w_i \wedge w_{i+2} \ldots w_k \\
&= (-1)^{gr(w_i)gr(w_{i+1})} \tau_k(a_1 \otimes w_1, \ldots \\
& \hspace{3.8 cm} \ldots, a_{i-1} \otimes w_{i-1}, a_{i+1} \otimes w_{i+1}, a_i \otimes w_i, a_{i+2} \otimes w_{i+2}, \ldots a_k \otimes w_k).
\end{split}
\end{equation*}

Além disso, dado $a \in A^\infty$ e $B = b \otimes w \in \Omega$, vale $\tau_1(aB) = \tau_1(Ba)$, já que 
\begin{equation*}
\begin{split}
\tau_1(aB) &= \tau_1(a (b\otimes w)) = \tau_1((a b)\otimes w) = \tau(ab) w \\
&= \tau(ba)w = \tau_1(ba \otimes w) = \tau_1(b \otimes w a) = \tau_1 (Ba).
\end{split}
\end{equation*} 

\fimdem

Além de $\tau_k$ estar bem definido um fato bastante importante é o de que 
$$ \tau( \delta (a)) = 0 \ \ \ , \ \ \ \forall a \in A^\infty ,$$
pois dado $a \in A^\infty$ e $X \in \mathfrak{g}$ com $g'_0 = X$, pela continuidade e pela $G$-invariância de $\tau$ temos:
\begin{equation} \label{taudelta}
\begin{split}
\tau( \delta_X (a)) = \tau \left( \underset{t \rightarrow 0}{\lim} \frac{\alpha_{g_t}(a)- a}{t}\right) = \underset{t \rightarrow 0}{\lim} \frac{\tau(\alpha_{g_t} (a)) - \tau(a)}{t} = 0
\end{split}
\end{equation}

Desta igualdade decorre   
\begin{equation} \label{dtau}
\begin{split}
\tau_1 \circ d = d \circ \tau_1
\end{split}
\end{equation}
pois
 $$d(\tau_1(a \otimes w)) = d(\tau(a)w) = \tau(a)dw  $$
e
 $$  \tau_1(d(a \otimes w))  = \tau_1(\delta(a)w + a dw) = \tau(\delta(a))w + \tau(a) dw  = \tau(a)dw. $$

Outro fato que será utilizado na proposição que segue é o de que o caminho (homotopia) entre duas curvaturas dada pela curva 
\begin{displaymath}
\begin{array}{ccccc}
\Theta :   & [0,1]  & \rightarrow  & A^\infty \otimes \Lambda^2  &\\
& t & \mapsto  & \Theta_t \ \ = &\Theta_0 + p(d(t\Gamma) + t\Gamma \wedge t \Gamma)p \\ 
\end{array} 
\end{displaymath}
vide a Proposição \ref{curv}, é diferenciável em relação a $t$, observando que esta é a derivação usual de uma curva no espaço de Fréchet $A^\infty \otimes \Lambda^2$.
 
Além disso, $\frac{d}{dt}\Theta_t = 0 + p ( d(\Gamma) + \Gamma \wedge t \Gamma + t \Gamma \wedge \Gamma )p$
e $\tau_1$ é um traço graduado, como vimos no Lema \ref{traco}, logo
$$\tau_1 (\frac{d}{dt}\Theta_t) = \tau_1( p d \Gamma p).$$

Como $\Gamma \in p M_n(A^\infty) p$, então $d(\tau_1(\Gamma)) = d (\tau_1(p \Gamma p))$. Assim por (\ref{dtau}) temos
$$d(\tau_1(\Gamma)) = d (\tau_1(p \Gamma p)) = \tau_1(d(p  \Gamma p)) = \tau_1 (dp \Gamma p  + p d \Gamma p - p \Gamma dp).$$
Utilizando o fato de $pdpp=0$, visto no Lema \ref{identidades}, o final da demonstração do Lema \ref{traco} nos leva a concluir que
\begin{equation*}
\begin{split}
d(\tau_1(\Gamma)) &= \tau_1 (p d \Gamma p)  +  \tau_1(dp \Gamma p)  -  \tau_1(p \Gamma dp) = \tau_1 (p d \Gamma p) + \tau_1(p dp \Gamma)    -  \tau_1(\Gamma dp p) \\
    &= \tau_1 (p d \Gamma p) + \tau_1(p dp p \Gamma)    - \tau_1(\Gamma p dp p) = \tau_1(p d \Gamma p).
\end{split}
\end{equation*}

Portanto $$\tau_1 (\frac{d}{dt}\Theta_t) = d (\tau_1(\Gamma)).$$

Com as informações que temos até aqui podemos enfim enunciar e demonstrar o teorema a seguir que é essencial na definição do caráter de Chern-Connes.

\begin{Prop} \label{Homol}
Com as notações definidas até o momento podemos afirmar que a forma diferencial $\tau_{k} (\Theta, \Theta, \ldots, \Theta) \in \Lambda^{2k}\mathfrak{g}^*$ é fechada e sua classe de cohomologia independe da conexão escolhida em $M^\infty$.
\end{Prop}
\Dem
Como provaremos, a seguir, que a forma $\tau_{k} (\Theta, \Theta, \ldots, \Theta)$ independe da conexão tomada, basta verificarmos, a princípio, que $\tau_{k} (\Theta_0, \Theta_0, \ldots, \Theta_0)$ é uma forma fechada, ou seja, $d(\tau_{k} (\Theta_0, \Theta_0, \ldots, \Theta_0)) = 0$.

Da equação \ref{pdp} e das observações feitas após a definição de $\tau_k$, respectivamente, temos 
\begin{equation*}
\begin{split}
\tau_{k} (\Theta_0, \Theta_0, \ldots, \Theta_0) &= \tau_{k} (pdp \wedge dp, \ldots, p dp \wedge dp) = \tau_{} (p^k)dp \wedge dp \wedge \ldots \wedge dp \wedge dp \\
   &= \tau(p)(dp \wedge dp)^{k} = \tau_1(p (dp \wedge dp)^k).
\end{split}
\end{equation*}
Vimos também, equação (\ref{dtau}), que $\tau_1 \circ d = d \circ \tau_1$, então
\begin{equation*}
\begin{split}
d(\tau_{k} (\Theta_0, \Theta_0, \ldots, \Theta_0)) &= d(\tau_1(p dp^{2k}))= \tau_1(d( p dp^{2k})) \\ &= \tau_1(dp \wedge dp^{2k}) = \tau_1(dp^{2k+1}) = 0, 
\end{split}
\end{equation*}
isto porque pelo final da demonstração do lema \ref{traco} temos
\begin{equation*}
\begin{split}
\tau_{1}(dp^{2k+1}) &=  \tau_{1}((2p-1)^2 dp^{2k+1}) \\
                       &= \tau_{1}((2p-1)dp^{2k+1}(2p-1)) 
\end{split}
\end{equation*}
e pelo Lema $\ref{identidades}$ podemos concluir que  
\begin{equation*}
\begin{split}  
 \tau_{1}(dp^{2k+1})  &= - \tau_{1}(-(2p-1)dp^{2k+1}(2p-1)) \\
			     &= - \tau_1(-(-1)^{2k+1}(2p-1)^2 dp^{2k+1}) \\
                       &= -\tau_{1}(dp^{2k+1}).
\end{split}
\end{equation*}

Resta portanto provar que $\tau_{k} (\Theta, \Theta, \ldots, \Theta)$ independe da conexão tomada, ou o que tem o mesmo efeito e que será demonstrado a seguir, que 
$$\tau_{k} (\Theta, \Theta, \ldots, \Theta)- \tau_{k} (\Theta_0, \Theta_0, \ldots, \Theta_0)$$ é uma forma exata.

Lembrando que dada uma conexão $\nabla$ qualquer temos $\nabla = \nabla^0 + \Gamma$, como visto no lema $\ref{conexao}$, assim tomemos $\nabla_t = \nabla^0 + t \Gamma$, para $t \in [0,1]$, um caminho entre $\nabla$ e $\nabla^0$ cuja curvatura será simbolizada por $\Theta_t$. 

Utilizando $\Theta_t =  \Theta_0 + p (d(t \Gamma) + t \Gamma \wedge t \Gamma) p \in \Omega^2 $ como visto na Proposição \ref{curv} e também as observações feitas antes do enunciado deste teorema temos
\begin{equation*}
\begin{split}  
\frac{d}{dt} \tau_k(\Theta_t, \ldots,  \Theta_t) &= \underset{i=0}{\overset{k-1}{\sum}} \tau_k(\underbrace{\Theta_t, \ldots, \Theta_t}_i, \frac{d}{dt}\Theta_t,\underbrace{\Theta_t, \ldots, \Theta_t}_{k-i-1}) = k \tau_k(\frac{d}{dt}\Theta_t, \Theta_t, \ldots, \Theta_t) \\
& = k \tau_k(pd(\Gamma) p, \Theta_t, \ldots, \Theta_t) = d (k \tau_k(\Gamma, \Theta_t, \ldots, \Theta_t)). 
\end{split}
\end{equation*}

Como $\tau_{k} (\Theta, \Theta, \ldots, \Theta)- \tau_{k} (\Theta_0, \Theta_0, \ldots, \Theta_0)= \int_0^1 \frac{d}{dt} \tau_{k}(\Theta_t, \ldots, \Theta_t) \ dt$, como esta integração é um limite de somas e portanto comuta com $d$ e também do fato de $\frac{d}{dt} \tau_k(\Theta_t, \ldots,  \Theta_t) = d(k \tau_k(\Gamma, \Theta_t, \ldots, \Theta_t) )$, podemos concluir que esta forma fechada independe da conexão tomada em $M^\infty$.
\fimdem

Cabe observar ainda que dadas duas projeções $p$ e $q$ em $\mathcal{P}_\infty (A^\infty)$ equivalentes, isto é, $[p]_0 = [q]_0 \in K_0(A)$ temos $[\tau_k(\Theta_0, \ldots, \Theta_0)] = [ \tau_k(\theta_0, \ldots, \theta_0)]$, sendo $\Theta_0 $ e $\theta_0$ as curvaturas das conexões de Levi-Civita associadas as projeções $p$ e $q$, respectivamente, e $[ \ . \ ]$ denotando a classe de cohomologia das formas invariantes à esquerda do grupo de Lie $G$.

Note que $\tau \left( \left( \begin{array}{ll} 
a & 0 \\
0 & 0 \\
\end{array} \right) \right) = \tau (a)$, para qualquer $a \in A$, já que por abuso de notação 
$\tau \left( \left( \begin{array}{ll} 
a & 0 \\
0 & 0 \\
\end{array} \right) \right) = \tau \circ Tr \left( \left( \begin{array}{ll} 
a & 0 \\
0 & 0 \\
\end{array} \right) \right)$, onde $Tr$ é o traço usual das matrizes.

Se $p \sim_0 q$, por $5.2.12$ de \cite{WO}, sabemos que $\left( \begin{array}{ll} p & 0 \\ 0 & 0 \\\end{array} \right) \sim_h \left( \begin{array}{ll} q & 0 \\ 0 & 0 \\\end{array} \right)$, onde $\sim_h$ significa que estas matrizes são homotópicas por caminho. 
Portanto podemos tomar para cada $t \in [0,1]$ a curvatura $\Theta_t = p_t dp_t \wedge dp_t$  com  $\Theta_0 $ a própria $\Theta_0$ e $\Theta_1 = \theta_0$. 

Como visto nas observações feitas antes da Proposição \ref{Homol} sabemos que existe $\frac{d}{dt}\theta_t$, que denotaremos apenas por $\theta'$. Para $\theta_t = p_t dp_t \wedge dp_t$; onde $p_t$ é uma homotopia que, como no Lema 4 de \cite{Fe},   podemos supor derivável em relação a variável $t$; temos
$$ \theta' = p' dp \wedge dp p + p(dp' \wedge dp + dp \wedge dp') + pdp \wedge dp p' = p' \theta + \theta p' + p d(p' dp - dp p')p.$$
Note que estamos omitindo o índice $t$ para simplificar a notação e também que estamos denotando por $p'$ a derivação de $p_t$ em relação a $t$.  

Assim como para a Proposição \ref{Homol}, podemos utilizar o Lema $\ref{identidades}$ para provar que $pp'p = 0$, logo $\theta p' = p' \theta = 0$, e ainda $p$ comuta com $p'dp$ e $dpp'$ e portanto $pd(p'dp-dpp')p = pd(p(p'dp-dpp'))p$. Assim
$$ \theta' =  p' \theta + \theta p' + p d(p' dp - dp p')p = pd(p'dp - dpp')p.$$

Utilizando agora o fato de que $pd(\theta)p = 0$, análogo ao visto após o Lema $\ref{identidades}$, e também que $\tau_1 (\frac{d}{dt} (\theta^k)) = k \tau_1(\theta' \theta^{k-1})$, assim como utilizado no final da demonstração anterior, temos 
\begin{equation*}
\begin{split}
\frac{d}{dt} \tau_{k}(\Theta_t, \ldots, \Theta_t) &= k \tau_{1}(\theta' \theta^{k-1}) = 
k \tau_{1}(pd(p'dp - dp p')p\theta^{k-1}) \\
&= k  \tau_{1}(pd((p'dp-dpp')\theta^{k-1})p) = d(k \tau_{1}(p(p'dp-dpp')\theta^{k-1})).\\
\end{split}
\end{equation*}

Como $\tau_{k} (\Theta, \Theta, \ldots, \Theta)- \tau_{k} (\Theta_0, \Theta_0, \ldots, \Theta_0)= \int_0^1 \frac{d}{dt} \tau_{k}(\Theta_t, \ldots, \Theta_t) \ dt$ e $$\frac{d}{dt} \tau_{k}
(\Theta_t, \ldots, \Theta_t)= d(k \tau_{1}(p(p'dp-dpp')\theta^{k-1})),$$ então $\tau_{k} (\Theta,
 \Theta, \ldots, \Theta)- \tau_{k} (\Theta_0, \Theta_0, \ldots, \Theta_0)$ é exata e portanto 
$$[\tau_k(\Theta_0, \ldots, \Theta_0)] = [\tau_k(\theta_0, \ldots, \theta_0)]$$
como queríamos.

A Proposição \ref{Homol} e a observação acima nos permitem definir, como feito na geometria diferencial, o caráter de Chern 
$$Ch_\tau ([p]_0) = \left[ \underset{k=0}{\overset{\infty}{\sum}} \left(\frac{1}{2 \pi i}\right)^k \frac{1}{k!} \tau_{k}( \Theta, \ldots, \Theta) \right], $$
que pelo Lema $\ref{identidades}$ e pela proposição anterior também pode ser definido como
\begin{equation*}
\begin{split}
Ch_\tau ([p]_0) &= \left[ \underset{k=0}{\overset{\infty}{\sum}} \left(\frac{1}{2 \pi i}\right)^k \frac{1}{k!} \tau_{k}( \Theta_0,\ldots, \Theta_0) \right] \\
&= \left[ \underset{k=0}{\overset{\infty}{\sum}} \left(\frac{1}{2 \pi i}\right)^k \frac{1}{k!} \tau_{k}( pdp \wedge dp,\ldots, p dp \wedge dp) \right] \\
&= \left[ \underset{k=0}{\overset{\infty}{\sum}} \left(\frac{1}{2 \pi i}\right)^k \frac{1}{k!} \tau_{1}( p( dp \wedge dp)^k) \right] \in H_{\mathbb{C}}^*(\mathfrak{g}).
\end{split}
\end{equation*}
Note que $G$ e $\mathfrak{g}$ têm dimensão finita, logo $H_{\mathbb{C}}^* (\mathfrak{g}) = 0$, para $*$ maior que a dimensão de $\mathfrak{g}$ e portanto qualquer das somas acima é finita.

Podemos melhorar um pouco mais esta definição lembrando que o anel de cohomologia de uma álgebra de Lie é isomorfo ao anel de cohomologia das formas invariantes de seu grupo de Lie, isto é, $H^*(\mathfrak{g}) \simeq H^*(G)$, vide $10.1.6$ de \cite{LO}.
Além disso, a soma acima é sempre finita logo $Ch_\tau([p]_0) \in H^{\textrm{par}}(G)$.

Provemos agora que em cada classe de equivalência da equação acima o representante obtido pela projeção $p$ é na verdade uma forma real.

\begin{Lem} 
Dada $p \in M_n(A^\infty)$ projeção, então $\frac{1}{(2 \pi i)^k} \tau_{k} (p (dp \wedge dp)^k)$ é um forma real, ou seja, $Ch_\tau([p]_0) \in H_{\mathbb{R}}^{\textrm{par}}(G)$.
\end{Lem}
\Dem
Em primeiro lugar vamos mostrar que dado $da \wedge db \in A^\infty \otimes \Omega^2$, temos $(da \wedge db)^* = -(db^* \wedge da^*)$. 

De fato, dados $X, Y \in \mathfrak{g}$ e lembrando por \ref{representacao} que $\delta$ é uma representação, temos 
\begin{equation*}
\begin{split}
(da \wedge db)^*(X,Y) &= ((da \wedge db)(X,Y))^* = (\delta_X(a)\delta_Y(b)-\delta_Y(a)\delta_X(b))^* \\
&= \delta_Y(b^*)\delta_X(a^*)-\delta_X(b^*)\delta_Y(a^*) = - (\delta_X(b^*)\delta_Y(a^*)-\delta_Y(b^*)\delta_X(a^*)) \\
&= - (db^* \wedge da^*)(X,Y).
\end{split}
\end{equation*}

Portanto 
$$\overline{\tau_{1}(p (dp \wedge dp)^k)} = \tau_{1}((p(dp\wedge dp)^k)^*) = \tau_{1}((-1)^k(dp \wedge dp)^k p) = (-1)^k \tau_{1} (p(dp \wedge dp)^k),$$
 ou seja, $\frac{1}{(2 \pi i)^k}\tau_{1}(p (dp \wedge dp)^k) \in \mathbb{R}$.
\fimdem



Com isto podemos enfim definir o caráter de Chern-Connes para o caso par.
\begin{Def} \label{Chernpar}
Com as notações que apresentamos até o momento o caráter de Chern-Connes é o homomorfismo 
$Ch_\tau :  K_0(A) \rightarrow H^{\textrm{par}}_\mathbb{R} (G)$, dado por
\begin{equation*}
\begin{split}
Ch_\tau ([p]_0) &= \left[ \underset{k=0}{\overset{\infty}{\sum}} \left(\frac{1}{2 \pi i}\right)^k \frac{1}{k!} \tau_k( \Theta, \ldots, \Theta) \right] \\
&= \left[ \underset{k=0}{\overset{\infty}{\sum}} \left(\frac{1}{2 \pi i}\right)^k \frac{1}{k!} \tau_k( \Theta_0, \ldots, \Theta_0) \right] \\
&= \left[ \underset{k=0}{\overset{\infty}{\sum}} \left(\frac{1}{2 \pi i}\right)^k \frac{1}{k!} \tau_1( p (dp \wedge dp)^k) \right]. 
\end{split}
\end{equation*}
\end{Def}

É conveniente lembrar que as somas da definição acima são finitas, já que se anulam para os valores em que $2k$ é maior do que a dimensão do grupo $G$ e esta é finita.

\section{A extensão para $K_1(A)$: o caso ímpar}

Nosso objetivo nesta seção é o de mostrar que se na definição do caráter de Chern-Connes, dada na seção anterior, trocarmos o C$^*$-sistema dinâmico $(A, G, \alpha)$ pelo também C$^*$-sistema dinâmico $(A \otimes C(S^1), G \times S^1, \alpha')$, com a ação $\alpha'$ dada por $\alpha'_{(g,h)} (x \otimes f) = \alpha_g (x) \otimes f_h$ onde $f_h(t) = f(t-h)$, estenderemos $Ch_\tau$ para $ch_\tau : K_0(A) \oplus K_1(A) \rightarrow H^*_\mathbb{R} (G)$.

\begin{Obs} \label{KKK}
Como vimos em \ref{cs1a} um resultado clássico de K-teoria, que também pode ser visto com mais detalhes em $8.B$ de \cite{WO} e $9.4.1$ de \cite{BB}, é o de que a partir da seqüência exata cindida 
\begin{equation*} 
 0 \longrightarrow SA \longrightarrow C(S^1, A) \overset{}{\longrightarrow} A  \longrightarrow 0
\end{equation*}
obtemos $$K_0(C(S^1,A)) = K_0(A) \oplus K_0(SA) \simeq K_0(A) \oplus K_1(A)$$ e como $C(S^1) \otimes A \simeq C(S^1,A)$, temos 
$$K_0(C(S^1) \otimes A) \simeq K_0(A) \oplus K_1(A).$$
\end{Obs}

\begin{Obs} \label{HHH}
Da mesma forma, pela fórmula de Künneth para o produto de cohomologias, vide exemplo $3B.3$ de \cite{H1}, temos $H^k (G \times S^1) \simeq H^k(G) \oplus H^{k-1}(G)$, isomorfismo este dado por
$$\begin{array}{ccc} H^k(G) \oplus H^{k-1}(G) & \rightarrow & H^k (G \times S^1)  \\
                    (a,b) & \mapsto & \pi_1^*(a)+\pi_1^*(b)\wedge d \theta \\
\end{array},$$
ressaltando que $\pi_1^*$ é o \emph{pull-back} da projecção $\pi_1 : G \times S^1 \rightarrow G$.

Com isto, temos $H^{2k}(G \times S^1) \simeq H^{2k}(G) \oplus H^{2k-1}(G)$ e como $H^{k}(G) = 0$, para $k < 0$, então 
$$H^{par} (G \times S^1) \simeq H^{*}(G).$$
\end{Obs}

Levando-se em conta os isomorfismos dados pela Observação \ref{KKK} e pela Observação \ref{HHH}, ao efetuarmos a mudança do C$^*$-sistema dinâmico $(A, G, \alpha)$ pelo $(A \otimes C(S^1), G \times S^1, \alpha')$ na Definição \ref{Chernpar}, teremos 
$$ch_\tau: K_0(A) \oplus K_1(A) \simeq K_0(A \otimes C(S^1)) \rightarrow  H^{\textrm{par}}(G\times S^1)\simeq H^*(G).$$

Note que $ch_\tau$ restrito a $K_0(A)$ é exatamente $Ch_\tau: K_0(A) \rightarrow H^{\textrm{par}}(G)$. Por este motivo iremos daqui para a frente nos preocupar com o que acontece no caso ímpar, isto é, com o que acontece com $ch_\tau$ quando calculado nas classes de unitários $[u]_1$ de $K_1(A)$.





Seja $[u]_1$ uma classe de equivalência em $K_1(A)$ e portanto $u \in \mathcal{U}_n(\mathcal{A})$, para algum $n \in \mathbb{N}$. Utilizando o isomorfismo $\theta_\mathcal{A}$ entre $K_1(A)$ e $K_0(SA)$, como pode ser visto em \ref{thetaa} ou também em $7.2.5$ de \cite{WO}, sabemos que $[u]_1 \in K_1(A), $ de dimensão $n \in \mathbb{N}$, é levado em $[p]_0-[p_n]_0$, com $p$ e $p_n$ projeções em $\mathcal{P}_{2n}(SA)$ \footnote{Estamos utilizando aqui $SA $ visto como $\{ f \in C([0, 2 \pi], A) \ : \ f(0) = f(2 \pi) = 0\}$.} dadas por 
$$p_n  = \left(  \begin{array}{ll} 
1_n & 0 \\
0 & 0_n \\
\end{array} \right)$$
e $p:  [0, 2 \pi] \ni \theta \mapsto p_\theta$,  com $p_\theta = w_\theta p_n w_\theta^*$ para 
$$ w_\theta = \left(  \begin{array}{ll} 
u & 0 \\
0 & 1_n \\
\end{array} \right)  \left(  \begin{array}{ll} cos(\frac{\theta}{4}) & -sen(\frac{\theta}{4}) \\
sen(\frac{\theta}{4}) & cos(\frac{\theta}{4}) \\
\end{array} \right)\left(  \begin{array}{ll} 
u^* & 0 \\
0 & 1_n \\
\end{array} \right)  \left(  \begin{array}{ll} cos(\frac{\theta}{4}) & sen(\frac{ \theta}{4}) \\
-sen(\frac{\theta}{4}) & cos(\frac{\theta}{4}) \\
\end{array} \right),$$ 
efetuando-se as multiplicações matriciais e utilizando algumas simplificações trigonométricas, temos em suma
$$p_\theta = \left(  \begin{array}{ll} 
1_n - \frac{sen^2( \theta/2)}{4} (1-u)(1-u^*) & (u-1)\frac{sen(\theta/2)}{2}\left[cos^2\left(\frac{ \theta}{4}\right) + u sen^2\left(\frac{\theta}{4}\right)\right] \\
(u^*-1)\frac{sen(\theta/2)}{2}\left[cos^2\left(\frac{\theta}{4}\right) + u^* sen^2\left(\frac{ \theta}{4}\right)\right] & \frac{sen^2(\theta/2)}{4} (1-u)(1-u^*) \\
\end{array} \right).$$

Assim através destes isomorfismos, ou melhor, através da Observação \ref{KKK} e da Observação \ref{HHH} temos para $u$, $p$ e $p_n$, como acima,
\begin{equation*} 
\begin{split}
ch_\tau: K_1(\mathcal{A}) & \hookrightarrow K_0(C(S^1)\otimes \mathcal{A}) \rightarrow H^{\textrm{par}}(G\times S^1) \simeq H^*(G) \\
[u]_1 & \mapsto  [p]_0 - [p_n]_0 \mapsto Ch_\tau([p]_0 - [p_n]_0).  
\end{split}
\end{equation*}
E portanto para $[u]_1 \in K_1(\mathcal{A})$ temos
$$ch_\tau([u_1]) = Ch_\tau([p]_0 - [p_n]_0) \in H^{*}(G).$$

\ \ \

\newpage

\ \ \ 

\newpage

\chapter{O caráter de Chern-Connes de $\overline{\Psi_{cl}^0(S^1)}$}
Neste capítulo explicitaremos o caráter de Chern-Connes, mais precisamente a matriz da transformação dada por este homomorfismo, definido em \cite{AC1} e apresentado em detalhes no capítulo anterior, no caso em que o C$^*$-sistema dinâmico em questão é $(\overline{\Psi_{cl}^0(S^1)} , S^1, \alpha)$, onde a ação $\alpha$ é a de conjugação pela translação no grupo de Lie $S^1$ e $\overline{\Psi_{cl}^0(S^1)}$ é a C$^*$-álgebra gerada pelos operadores pseudodiferenciais clássicos de ordem zero na variedade $S^1$.

Cabe ressaltar que a ação $\alpha$ tomada é de classe $C^\infty$ em $\Psi_{cl}^0(S^1)$ e portanto contínua em seu fecho, isto é, contínua em $\overline{\Psi_{cl}^0(S^1)}$. Além disso, $\Psi_{cl}^0(S^1)$ é invariante pelo cálculo funcional holomorfo e portanto podemos tomar a álgebra $A^\infty$ utilizada na definição do caráter de Chern-Connes como sendo $\Psi_{cl}^0(S^1)$, ou seja, estamos de fato calculando o caráter de Chern-Connes para a álgebra dos operadores pseudo\-diferenciais clássicos de ordem zero de $S^1$.

Como visto no capítulo anterior, para calcularmos o homomorfismo de Chern-Connes precisamos antes conhecer os grupos de K-teoria da C$^*$-álgebra $\overline{\Psi_{cl}^0(S^1)}$ e seus geradores e para isto utilizaremos alguns resultados clássicos de  "Comparison Algebras" \ que podem ser vistos em VI.2 de \cite{CO} ou em \cite{TO} e mais especificamente, neste caso em que a variedade é $S^1$, pode-se consultar também \cite{TO0}. 

Um destes resultados é o de que $\overline{\Psi_{cl}^0(S^1)}$ é isomorfa a C$^*$-álgebra $\mathcal{A}$ que definiremos a seguir. Mas para tal lembremos que $a(D_\theta) = F_d^{-1}M_aF_d$, onde $M_a$ é o operador de multiplicação pela seqüência $(a_i) \in CS(\mathbb{Z})$, $F_d$ designa a transformada de Fourier discreta, isto é, $F_d: L^2 (S^1) \rightarrow l^2(\mathbb{Z})$, onde $(F_du)_j = \frac{1}{\sqrt{2 \pi}} \int_{- \pi}^\pi \tilde{u}(\theta) e^{-ij \theta} d \theta$, $j \in \mathbb{Z}$ e $\tilde{u}(\theta) = u (e^{i \theta})$ e $CS(\mathbb{Z}) = \{ (a_i)_{i \in \mathbb{Z}} \ : \ \underset{i \rightarrow \infty}{\lim} a_i = a(\infty) \ \ \textrm{e} \ \  \underset{i \rightarrow  -\infty}{\lim} a_i = a(-\infty) \ \  \textrm{existem} \}$. 

\begin{Def} \label{SS1}
Denotaremos por $\mathcal{A}$ a C$^*$-subalgebra de $\mathcal{L}(L^2(S^1))$, operadores li\- mitados de $L^2(S^1)$, gerada por $\mathcal{A}_1$ e $\mathcal{A}_2$, com $\mathcal{A}_1 = \{ M_a : a \in C^\infty(S^1) \}$ e $\mathcal{A}_2=\{ b(D_\theta) : b \in CS(\mathbb{Z})\}$.
\end{Def}

Um fato bastante conhecido e que pode ser visto em \cite{TO} é o de que $\mathcal{K} \subset \mathcal{A}$ e que a sequência
\begin{equation} \label{seqS1}
 0 \longrightarrow \mathcal{K} \longrightarrow \mathcal{A} \overset{\pi}{\longrightarrow} \mathcal{A}/\mathcal{K}  \longrightarrow 0
\end{equation}
é exata. E ainda, pelo Teorema 2 de \cite{TO} e seu corolário, sabemos que
existe um isomorfismo $\varphi$ entre o quociente $\mathcal{A}/\mathcal{K}$ e $C(S^*S^1) \overset{\psi}{\simeq} C(S^1 \times \{-\infty, +\infty \})$. Como dito anteriormente isto é uma conseqüência de um resultado mais geral sobre "Compa\-rison Algebras" \ que pode ser visto em VI.2 de \cite{CO}.

Ao compormos $\varphi$ com $\pi$, obtemos o isomorfismo $\sigma = \varphi \circ \pi $ conhecido como símbolo principal. E finalmente, compondo $\sigma$ com $\psi$ obtemos um isomorfismo que transforma os geradores da seguinte maneira, 

\begin{equation} \label{gerador1}
\psi \circ \sigma([M_a]) = \psi([\sigma_{M_a}]_\mathcal{K}) = (a , a) \ \ , \ \textrm{ para } \ a \in C^\infty (S^1) \ \ \ \textrm{ e }
\end{equation}

\begin{equation} \label{gerador2}
\psi \circ \sigma([b(D_\theta)])= \psi([\sigma_{b(D_\theta)}]_\mathcal{K}) = (b(-\infty), b(\infty)) \ \ \textrm{,    para } (b_j) \in CS(\mathbb{Z}).
\end{equation}

Pela periodicidade de Bott, vide \ref{Bott}, a partir de uma seqüência exata curta obtemos a seqüência exata de seis termos da K-teoria complexa, como pode ser visto em \ref{seis}. E em nosso caso, isto é, para a seqüência (\ref{seqS1}), obtemos a seqüência exata cíclica dada por
$$ \xymatrix{
K_0(\mathcal{K})  \ar[rr]^-{}  & \ \ &  K_0(\mathcal{A}) \ \ \ar[rr]^-{}     & \ \ &   K_0(\mathcal{A} / \mathcal{K}) \ar[dd]^-{\delta_0}     \\
& \\
K_1(\mathcal{A} / \mathcal{K}) \ar[uu]^-{\delta_1} & \ \ &  K_1(\mathcal{A})  \ar[ll]_-{}  & \ \  & K_1(\mathcal{K}) \ar[ll]_-{}                     }.$$
\begin{equation}\label{seis1}\end{equation}

Como pode ser visto em \ref{compacto} ou também na tabela de grupos de K-teoria que consta ao final de \cite{RLL} temos
$$K_0(\mathcal{K}) = \mathbb{Z} \ \ \ \ \ \textrm{ e } \ \ \ \ \ K_1(\mathcal{K})= 0$$
e como $\mathcal{A}/\mathcal{K} \simeq C(S^1 \times \{-\infty, +\infty \})$, por \ref{Ks1} temos portanto
$$K_0(\mathcal{A}/\mathcal{K}) \simeq  K_0(C(S^1 \times \{-\infty, +\infty \})) = \mathbb{Z} \oplus \mathbb{Z} \ \ \ \textrm{ e }$$ 
$$K_1(\mathcal{A}/\mathcal{K}) \simeq  K_1(C(S^1 \times \{-\infty, +\infty \})) = \mathbb{Z} \oplus \mathbb{Z}.$$

Quando $\mathcal{K}$ aparece na seqüência exata curta, como em nosso exemplo \ref{seqS1}, a aplicação $\delta_1: K_1(\mathcal{A}/\mathcal{K})\rightarrow K_0(\mathcal{K}) = \mathbb{Z}$ é o índice de Fredholm, vide $8.3.2$ de \cite{BB} ou os comentários após \ref{indice}. E, além disso, neste caso particular $\delta_1$ é sobrejetora, como mostra o seguinte lema.

\begin{Lem} \label{Gerador}
A aplicação do índice $\delta_1 : K_1(\mathcal{A}/\mathcal{K}) \rightarrow K_0(\mathcal{K})$ é sobrejetora. 
\end{Lem}
\begin{Dem}
Tomemos o operador $B \in \mathcal{A}$, dado por 
$$Bu(z) = \underset{j \in \mathbb{Z}}{\sum} z^j b_j (z) \hat u_j , \ \ \ \textrm{onde} \ \ b_j(z) = \left\{ \begin{array}{ll}
1  &, \textrm{se $j \geqslant 0$} \\
z  &, \textrm{se $j<0$} \\
\end{array} \right..$$ 

Lembrando que $\cal{A}$ é gerada por $\mathcal{A}_1$ e $\mathcal{A}_2$, como na Definição \ref{SS1}, podemos mostrar que o operador $B$ está de fato em $\cal{A}$, pois $B = M_a(I-H(D_\theta))+H(D_\theta)$ para 
$a(z) = z$ e $H(j) = \left\{ \begin{array}{ll}
1  &, \textrm{se $j \geqslant 0$} \\
0  &, \textrm{se $j<0$} \\
\end{array} \right..$  

Assim, $F_d B F_d ^{-1}(u_j)_{j \in \mathbb{Z}}\mid_k = \left\{ \begin{array}{lll}
u_{k-1}  &, \textrm{se $k < 0$} &\\
u_0+u_{-1}  &, \textrm{se $k=0$} &\\
u_{k}  &, \textrm{se $k > 0$} &\\
\end{array} \right.$, ou seja, $$\delta_1 ([[B]_\mathcal{K}]_1) = ind([B]_\mathcal{K}) = dim(ker B) - codim(Im B) = 1 - 0 = 1$$ e portanto existe um operador que é levado no gerador do grupo $K_0(\mathcal{K})$, logo a aplicação $\delta_1 $ é sobrejetora.
 \fimdem \end{Dem}

Com base nestas informações, isto é, nas equações (\ref{gerador1}) e (\ref{gerador2}) e no lema anterior, o diagrama (\ref{seis1}) pode ser entendido como

$$ \xymatrix{
\mathbb{Z}  \ar[rr]^-{0}  & \ \ &  K_0(\mathcal{A}) \ \ \ar[rr]^-{}     & \ \ &   \mathbb{Z} \oplus \mathbb{Z} \ar[dd]^-{\delta_0}     \\
& \\
\mathbb{Z} \oplus \mathbb{Z} \ar[uu]^-{\delta_1} & \ \ &  K_1(\mathcal{A})  \ar[ll]_-{}  & \ \  & 0 \ar[ll]_-{}                     }. $$

Denotaremos por $(f, g)$ a função de $C(S^1 \times \{-\infty, \infty \})$ em que $f \in C(S^1)$ está sobre $S^1 \times \{ -\infty \} $ e $g \in C(S^1)$ sobre $S^1 \times \{ \infty \}$ e também, como em \ref{Ks1}, designamos por $\mathfrak{z}$, $\mathfrak{0}$ e $\mathfrak{1}$, respectivamente, as funções dadas por $\mathfrak{z} (z) = z$, $\mathfrak{0} (z) = 0$ e $\mathfrak{1} (z) = 1$, para todo $z \in S^1$. 

Pela seqüência exata cíclica acima temos portanto que 
$$K_0(\mathcal{A}) \simeq K_0 (\mathcal{A}/\mathcal{K}) \simeq K_0(C(S^*S^1)) = [(\mathfrak{1},\mathfrak{0})]_0 \mathbb{Z} \oplus [(\mathfrak{0},\mathfrak{1})]_0 \mathbb{Z}.$$ 
Mais precisamente, pelo Lema \ref{Gerador} temos $\sigma_{H(D_\theta)} = (\mathfrak{1},\mathfrak{0})$, logo
$$ K_0(\mathcal{A}) = [[H(D_\theta)]_\mathcal{K}]_0 \mathbb{Z} \oplus [[I-H(D_\theta)]_\mathcal{K}]_0 \mathbb{Z}.$$ 

E também, como $\delta_1([(\mathfrak{1},\mathfrak{1})]) = 0$ e $K_1(\mathcal{A}) \simeq ker \delta_1 $, temos 
$$K_1(\mathcal{A}) = [[\mathfrak{z}]_\mathcal{K}]_1 \mathbb{Z},$$
 já que $\sigma_\mathfrak{z} = (\mathfrak{1}, \mathfrak{1})$. 

Utilizando os grupos de K-teoria que acabamos de calcular e o fato, bastante conhecido, de que $H_\mathbb{R}^* (S^1) = H_\mathbb{R}^0(S^1) \oplus H_\mathbb{R}^1(S^1) = \mathbb{R} \oplus \mathbb{R}$, o cárater de Chern-Connes, para nosso C$^*$-sistema dinâmico $(\overline{\Psi_{cl}^0(S^1)} \simeq \mathcal{A}, S^1, \alpha)$,  é dado particularmente por 
$$ch_\tau : [[H(D_\theta)]_\mathcal{K}]_0 \mathbb{Z} \oplus [[I-H(D_\theta)]_\mathcal{K}]_0 \mathbb{Z} \oplus [[\mathfrak{z}]_\mathcal{K}]_1 \mathbb{Z} \rightarrow \mathbb{R} \oplus \mathbb{R}.$$

Como para $k \geqslant 2$ temos $H^k(S^1) = 0$, neste nosso exemplo o caso par, isto é, em $K_0(\mathcal{A})$ o caráter de Chern-Connes se resume a parcela do somatório em que $k=0$. Portanto dada uma projeção $p$ de $\mathcal{A}$ a definição
$$Ch_\tau ([p]_0) = \underset{k=0}{\overset{\infty}{\sum}} \frac{1}{(2 \pi i)^k k!} \tau_k (p(dp)^{2k}),$$
se resume simplesmente a $Ch_\tau ([p]_0) = \underset{k=0}{\overset{0}{\sum}} \frac{1}{(2 \pi i)^k k!} \tau_k (p(dp)^{2k}) = \tau(p)$.

Já o caso ímpar, isto é, o caráter de Chern-Connes calculado em $K_1(\mathcal{A})$ também não se anula apenas para a parcela do somatório em que $k=1$. Por isto, dado um unitário $u$ de $\mathcal{A}$, vamos provar que $ch_\tau (u) = \frac{1}{2 \pi i} \tau( u^* \delta(u) )$ e posteriormente calcular este valor.

Na verdade o que vamos provar agora é um resultado mais geral que vale para qualquer grupo a um parâmetro, como afirmado e não demonstrado no \emph{caso parti\-cular} (\emph{a}) de \cite{AC1}, por isso não faremos referências explícitas ao nosso exemplo parti\-cular para provar que dado um unitário $u$ da álgebra temos $ch_\tau([u]_1) = \frac{1}{2 i \pi} \tau(u^* \delta(u)) $ desde que o grupo em questão seja a um parâmetro.

Na discussão feita na segunda seção do capítulo anterior vimos que 
$$ch_\tau ([u]_1) = Ch_\tau([p]_0 - [p_n]_0 ),$$ 
onde $Ch_\tau([p]_0-[p_n]_0)$ é o caráter de Chern-Connes como dado na Definição $\ref{Chernpar}$.

Antes de começarmos os cálculos faremos as seguintes observações.

\begin{Obs} \label{k02}
O primeiro termo do somatório de $Ch_\tau([p]_0-[p_n]_0)$ é nulo. De fato, para $k=0$ temos $\tau_0(p-p_n) = \tau(p - p_n)$ que por abuso de notação representa $\tau \circ Tr (p - p_n)$ e como 
$$p - p_n = \left(  \begin{array}{ll} 
- \frac{sen^2( \theta/2)}{4} (1-u)(1-u^*) & (u-1)\frac{sen(\theta/2)}{2}
\left[cos^2\left(\frac{ \theta}{4}\right) + u sen^2\left(\frac{\theta}{4}\right)\right] \\
(u^*-1)\frac{sen(\theta/2)}{2}\left[cos^2\left(\frac{\theta}{4}\right) + 
u^* sen^2\left(\frac{ \theta}{4}\right)\right] & \frac{sen^2(\theta/2)}{4} (1-u)(1-u^*) \\
\end{array} \right) $$
 temos portanto
$$\tau \circ Tr (p - p_n) = \tau (Tr(p - p_n)) = 0.$$

Além disso, para $k > 2$ sabemos que $H^k(G) = 0$ e como $dp_n = 0$ temos apenas 
$$ ch_\tau ([u]_1) = Ch_\tau([p]_0 - [p_n]_0 ) = 0 + \tau_1(p dp \wedge dp) + 0.$$
\end{Obs}

Para a observação que faremos abaixo estamos utilizando os seguintes fatos  
$$dp (X,0) = \delta_{(X,0)}(p): S^1 \ni \theta  \mapsto \delta_X(p_\theta)$$
e  
$$dp (0, d\theta) = \delta_{(0, d \theta)}(p) S^1 \ni \theta  \mapsto \frac{\partial p}{\partial \theta}.$$ 
Além disso, cabe ressaltar que $S^1$, aqui parametrizado por $S^1 = \{ e^{i \theta}: \theta \in \mathbb{R} \}$, não é nosso exemplo particular e sim $S^1$ da definição do caráter de Chern-Connes visto na segunda seção do capítulo anterior.

\begin{Obs} \label{elim}
Utilizando o isomorfismo $H^2(G \times S^1) \simeq H^2(G) \oplus H^1(G)$, como na Observação \ref{HHH}, e o fato de que $H^2(G) = 0$ temos $H^2(G \times S^1) \simeq H^1(G)$ e através disto podemos concluir que $\tau_1(pdp \wedge dp)$ apresenta o termo $d\theta$. Em outras palavras, de todas as  possibilidades para se calcular $dp \wedge dp$ em combinações do tipo $(X, d\theta) \in \mathfrak{g} \times \mathfrak{s}^1$, temos
$$ dp \wedge dp ((X, 0), (X,0)) = \delta_X(p)\delta_X(p) - \delta_X(p)\delta_X(p) = 0$$
e
$$ dp \wedge dp ((0, d\theta), (0 , d\theta)) = \frac{\partial p }{\partial \theta} \frac{\partial p }{\partial \theta} - \frac{\partial p }{\partial \theta} \frac{\partial p }{\partial \theta} = 0.$$

Assim  basta calcularmos $pdp \wedge dp ((X,0) , (0, d\theta))$, isto é,  
$$p dp \wedge dp ((X , 0), (0, d \theta)) = p \left( \delta_X(p) \frac{\partial p}{\partial \theta}  - \frac{\partial p}{\partial \theta}  \delta_X (p) \right).$$

\end{Obs}

Como 
$$p_\theta = \left(  \begin{array}{ll} 
1_n- \frac{sen^2( \theta/2)}{4} (1-u)(1-u^*) & (u-1)\frac{sen(\theta/2)}{2}\left[cos^2\left(\frac{ \theta}{4}\right) + u sen^2\left(\frac{\theta}{4}\right)\right] \\
(u^*-1)\frac{sen(\theta/2)}{2}\left[cos^2\left(\frac{\theta}{4}\right) + u^* sen^2\left(\frac{ \theta}{4}\right)\right] & \frac{sen^2(\theta/2)}{4} (1-u)(1-u^*) \\
\end{array} \right)$$
temos 
$\delta_X(p) = \left(  \begin{array}{ll} \delta_X(p)_{11} & \delta_X(p)_{12} \\ \delta_X(p)_{21} & \delta_X(p)_{22} \\ \end{array} \right)$,
com
\begin{itemize}
\item $\delta_X(p)_{11} = \frac{sen^2( \theta/2)}{4} (\delta(u)+ \delta(u^*) ) , $

\item $\delta_X(p)_{12} = \frac{sen(\theta/2)}{2}\left[cos^2\left(\frac{ \theta}{4}\right) \delta(u) + sen^2\left(\frac{\theta}{4}\right)(\delta(u) u + u \delta(u) )\right], $

\item $\delta_X(p)_{21} = \frac{sen(\theta/2)}{2}\left[cos^2\left(\frac{\theta}{4}\right) \delta(u^*) + sen^2\left(\frac{ \theta}{4}\right) (\delta(u^*)u^* + u^* \delta(u^*)) \right] $ e
 
\item $ \delta_X(p)_{22} = - \delta_X(p)_{11} = - \frac{sen^2(\theta/2)}{4} (\delta(u)+ \delta(u^*)).$
\end{itemize}

E $\frac{\partial p}{\partial \theta} = \left(  \begin{array}{ll} \frac{\partial p}{\partial \theta}_{11} & \frac{\partial p}{\partial \theta}_{12} \\ \frac{\partial p}{\partial \theta}_{21} & \frac{\partial p}{\partial \theta}_{22} \\ \end{array} \right)$,
onde
\begin{itemize}
\item $ \frac{\partial p}{\partial \theta}_{11} = - \frac{sen( \theta)}{8} (1-u)(1-u^*) , $

\item $ \frac{\partial p}{\partial \theta}_{22} = - \frac{\partial p}{\partial \theta}_{11} = \frac{sen( \theta)}{8} (1-u)(1-u^*) , $

\item $\frac{\partial p}{\partial \theta}_{12} = (u-1) \frac{cos(\theta/2)}{4}\left[cos^2\left(\frac{ \theta}{4}\right) + u sen^2\left(\frac{\theta}{4}\right)\right]  + (u-1) \frac{sen(\theta/2)}{2} \left[\frac{-1}{4}sen\left(\frac{ \theta}{2}\right) + u \frac{1}{4}sen\left(\frac{\theta}{2}\right)\right] $ e
\item $\frac{\partial p}{\partial \theta}_{21} = \frac{\partial p}{\partial \theta}_{12}^* .$
\end{itemize}

Sendo assim,    
$p dp \wedge dp = p \left( \delta_X(p) \frac{\partial p}{\partial \theta}  - \frac{\partial p}{\partial \theta}  \delta_X (p) \right)  = $ \\ 
$$p \left[  \left(  \begin{array}{ll} \delta_X(p)_{11} & \delta_X(p)_{12} \\ \delta_X(p)_{21} & \delta_X(p)_{22} \\ \end{array} \right) \left(  \begin{array}{ll} \frac{\partial p}{\partial \theta}_{11} & \frac{\partial p}{\partial \theta}_{12} \\ \frac{\partial p}{\partial \theta}_{21} & \frac{\partial p}{\partial \theta}_{22} \\ \end{array} \right) +  \left(  \begin{array}{ll} \frac{\partial p}{\partial \theta}_{11} & \frac{\partial p}{\partial \theta}_{12} \\ \frac{\partial p}{\partial \theta}_{21} & \frac{\partial p}{\partial \theta}_{22} \\ \end{array} \right) \left(  \begin{array}{ll} \delta_X(p)_{11} & \delta_X(p)_{12} \\ \delta_X(p)_{21} & \delta_X(p)_{22} \\ \end{array} \right) \right] $$ \\
após feitas as devidas multiplicações, que em nosso caso foram computadas com o auxílio de um programa computacional matemático e posteriormente conferidas \footnote{Como esses cálculos são muito extensos não os apresentamos neste texto, mas aqueles que desejarem vê-los podem consultá-los na página virtual cujo endereço é \emph{www.ime.usp.br/$\sim$dpdias} .}, tomando-se o traço $Tr(p dp \wedge dp)$ obtemos um somatório em que todos os termos, exceto $u^* \delta(u)$, podem ser escritos na forma $\delta (u^m)$, para algum $m$ inteiro. E como $\tau \circ \delta = 0$, vide as explicações dadas para a equação (\ref{taudelta}), temos portanto 
\begin{equation} \label{impar1}
\tau(p dp \wedge dp) = \tau(Tr(p dp \wedge dp)) = \tau(u^* \delta u).
\end{equation}

Voltando agora ao nosso C$^*$-sistema dinâmico $(\overline{\Psi_{cl}^0(S^1)},S^1, \alpha)$ tomemos o traço $\tau = \tau_\mathcal{K} \circ \psi \circ \sigma$, onde  as aplicações $\sigma$ e $\psi$ são as dadas pelas equações (\ref{gerador1}) e (\ref{gerador2}) e $\tau_\mathcal{K}$ é o traço de $C(S^1 \times \{-\infty, \infty\}) \simeq \mathcal{A}/\mathcal{K}$ dado pela combinação linear das integrais sobre cada uma das cópias de $S^1$, isto é, $S^1 \times \{- \infty \}$ e $S^1 \times \{ \infty \}$.

Simbolicamente temos
$$\tau (a) = \frac{c_1}{2 \pi}  \int_{S^1} \sigma_a(z, \infty) dz + \frac{c_2}{2 \pi} \int_{S^1} \sigma_a(z, -\infty) \ dz, \ \ a \in \mathcal{A}.$$ 

Este traço, assim definido, é um funcional linear contínuo e invariante pela ação de translação $\alpha$, mas para que $\tau(a^*) = \overline{\tau(a)}$, como na Definição \ref{tracog}, é necessário que $c_1$ e $c_2$ sejam números reais.

Só nos resta então calcular o traço dado em cada um dos geradores do grupo, a fim de obter a matriz da transformação dada pelo caráter de Chern-Connes. Portanto, para $c_1$ e $c_2$ em $\mathbb{R}$, temos
$$\tau([[H(D_\theta)]_\mathcal{K}]_0) = \frac{c_1}{2 \pi} \int_{S^1} \sigma_{H(D_\theta)}(z, \infty) dz + \frac{c_2}{2 \pi} \int_{S^1} \sigma_{H(D_\theta)}(z, -\infty) \ dz = c_1 .$$ 

De forma análoga $\tau([[I-h(D_\theta)]_\mathcal{K}]_0) = c_2 $ e 
$$\tau(\mathfrak{z}^{-1} d(\mathfrak{z}) ) = \frac{c_1}{2 \pi} \int_{S^1} \sigma_{\mathfrak{z}}(z, \infty) dz + \frac{c_2}{2 \pi} \int_{S^1} \sigma_{\mathfrak{z}}(z, -\infty) \ dz = 0.$$ 

Unindo-se as informações acima podemos afirmar que a matriz da transformação dada pelo homomorfismo de Chern-Connes para o C$^*$-sistema dinâmico $(\overline{\Psi_{cl}^0(S^1)} , S^1, \alpha)$  dado por
$$ch_\tau : [[H(D_\theta)]_\mathcal{K}]_0 \mathbb{Z} \oplus [[I-H(D_\theta)]_\mathcal{K}]_0 \mathbb{Z} \oplus [[\mathfrak{z}]_\mathcal{K}]_1 \mathbb{Z} \rightarrow \mathbb{R} \oplus \mathbb{R}$$
 é 
$$ \left( \begin{array}{ccc}
c_1  & c_2   & 0\\
0 & 0 & 0\\
\end{array} \right).$$
Isto pois
 $$ \tau([[H(D_\theta)]_\mathcal{K}]_0) = \left( \begin{array}{ccc}
c_1  & c_2   & 0\\
0 & 0 & 0\\
\end{array} \right)  \left( \begin{array}{lll}
1 \\
0 \\
0  \\
\end{array} \right) =  \left( \begin{array}{cc}
c_1 \\
0 \\
\end{array} \right),$$  
$$ \tau ([[I-H(D_\theta)]_\mathcal{K}]_0) = \left( \begin{array}{ccc}
c_1  & c_2   & 0\\
0 & 0 & 0\\
\end{array} \right)  \left( \begin{array}{lll}
0 \\
1 \\
0  \\
\end{array} \right) =  \left( \begin{array}{cc}
c_2 \\
0 \\
\end{array} \right)$$
 e 
$$\tau ([\mathfrak{z}]_1) = \left( \begin{array}{ccc}
c_1  & c_2   & 0\\
0 & 0 & 0\\
\end{array} \right)  \left( \begin{array}{lll}
0 \\
0 \\
1  \\
\end{array} \right) =  \left( \begin{array}{cc}
0 \\
0 \\
\end{array} \right).$$ 

\ \ \

\newpage

\ \ \ 

\newpage

\chapter{O caráter de Chern-Connes de $\overline{\Psi_{cl}^0(S^2)}$}
Nosso objetivo neste capítulo é o de calcular o caráter de Chern-Connes, assim como no capítulo anterior, porém agora para $(\overline{\Psi_{cl}^0(S^2)} , SO(3), \alpha)$, o C$^*$-sistema dinâmico onde $\overline{\Psi_{cl}^0(S^2)}$, que será denotada por $\mathcal{A}$, é a C$^*$-álgebra gerada pelos ope\-radores pseudodiferenciais clássicos de ordem zero da esfera e a ação $\alpha$ é a de conjugação pela translação no grupo de Lie $SO(3)$, em outra palavras, a ação $\alpha$ leva $g \in SO(3)$ em $T_gAT_g^{-1}$, com $T_g u(x) = u(g^{-1}x)$, para $u \in L^2(S^2)$. 

Novamente a ação $\alpha$ é contínua, pois é de classe $C^\infty$ em $\Psi_{cl}^0(S^2)$, como pode ser visto em \cite{CM}, e portanto contínua em $\overline{\Psi_{cl}^0(S^2)}$. E repetindo o raciocínio do início da capítulo dois temos $\Psi_{cl}^0(S^2)$ também invariante pelo cálculo funcional holomorfo 
, ou seja, podemos pensar em $\Psi_{cl}^0(S^2)$ como sendo a álgebra $A^\infty$ do primeiro capítulo e portanto podemos afirmar que estamos de fato calculando o caráter de Chern-Connes da álgebra dos operadores pseudodiferenciais clássicos de ordem zero da variedade $S^2$.

Para encontrarmos tal homomorfismo precisamos antes saber quais são os grupos de K-Teoria desta C$^*$-álgebra e assim como feito anteriormente, em (\ref{seqS1}), podemos utilizar a seqüência exata 
\begin{equation} \label{seqS2}
 0 \longrightarrow \mathcal{K} \longrightarrow \mathcal{A} \overset{\pi}{\longrightarrow} \mathcal{A}/\mathcal{K}  \longrightarrow 0
\end{equation}
e novamente, como conseqüência dos resultados de "Comparison Algebras" \ obtidos por Cordes em VI.2 de \cite{CO}, temos para este caso $\mathcal{A}/ \mathcal{K} \simeq C(S^*S^2) \simeq C(SS^2)$.

A fim de encontrar os grupos de K-teoria de $\mathcal{A}$, precisaremos saber de antemão quais são os grupos de K-teoria das funções contínuas sobre o fibrado das esferas ou das coesferas de $S^2$, já que estes são homeomorfos. A seção a seguir se dedica a estes cálculos, isto é, ao cálculo dos grupos de K-teoria de $C(SS^2)$ e também de seus geradores.

\section{A K-teoria do fibrado das esferas da esfera}

A seqüência de Mayer-Vietoris é uma ferramenta bastante conhecida e utilizada da Topologia Algébrica e este resultado pode ser estendido para qualquer Teoria de (co)homologia. A idéia básica, como pode ser visto na seção $7.2$ de \cite{SM}, é a de que a partir de um espaço $X$ Hausdorff e compacto que é a união de dois subespaços compactos que possuem uma certa inteseccção, como representado no diagrama 
$$ \xymatrix{
 X & X_1 \ar[l] \\
 X_2 \ar[u] & X_1 \cap X_2 \ar[u] \ar[l] } $$
onde as aplicações em questão são injetoras, podemos obter por dualidade o seguinte diagrama
$$ \xymatrix{
 C(X) \ar[d] \ar[r] & C(X_1) \ar[d] \\
 C(X_2) \ar[r] & C(X_1 \cap X_2)  } $$
onde as aplicações nada mais são do que restrições e com isto podemos mostrar que 
$$C(X) \simeq  \{ (a_1 , a_2) \in C(X_1) \oplus C(X_2) : a_{1 \mid_{X_1 \cap X_2}} = a_{2 \mid_{X_1 \cap X_2}}   \} .$$ 

Com base nestas idéias enunciamos o Teorema a seguir cuja demonstração pode ser encontrada em 7.2.1 de \cite{SM}.

\begin{Teo} \label{MVH}
Dado o diagrama comutativo de C$^*$-álgebras com unidade
$$\xymatrix{ A \ar[d]_{p_2} \ar[r]^{p_1}& B_1 \ar[d]^{\pi_1}\\
B_2 \ar[r]_{\pi_2}& D} $$
com $A$ o produto fibrado de $B_1$ e $B_2$ sobre $D$, isto é, 
$$A = \{ (b_1, b_2) : \pi_1(b_1) = \pi_2(b_2) \} \subseteq B_1 \oplus B_2 \ ,$$  
$\pi_1$ e $\pi_2$  $*$- homomorfismo sobrejetivos e $p_k$, $k=1 \textrm{ ou } 2$, as restrições das projeções $i_k : B_1 \oplus B_2 \rightarrow B_k$, temos a seqüência exata cíclica de seis termos
$$ \xymatrix{
K_0(A)  \ar[rr]^-{(p_{1_*},p_{2_*})}  & \ \ &  K_0(B_1) \oplus K_0 (B_2) \ \ \ar[rr]^-{\pi_{2_*} - \pi_{1_*}}     & \ \ &   K_0(D) \ar[dd]     \\
& \\
K_1(D) \ar[uu]^{}  & \ \ &  K_1(B_1) \oplus K_1(B_2) \ar[ll]_-{\pi_{2_*} - \pi_{1_*}}  & \ \  & K_1(A) \ar[ll]_-{(p_{1_*},p_{2_*})}                     } .$$

\end{Teo}

O teorema enunciado acima é um caso particular do Teorema 21.2.2 (Seqüência de Mayer-Vietoris para Homologia) de \cite{BB} sobre a existência de uma seqüência de Mayer-Vietoris para uma teoria de Homologia qualquer.



Nosso intuito é o de utilizar a seqüência de Mayer-Vietoris, vista no Teorema \ref{MVH}, em $C(S^*S^2)$, mas para facilitar os cálculos e a compreensão faremos uso desta seqüência em $C(SS^2)$, já que $C(SS^2) \simeq C(S^*S^2)$ como dito anteriormente. E para isto utilizaremos neste nosso caso particular, isto é, em $SS^2$ as observações feitas antes do Teorema \ref{MVH} como podemos ver no lema a seguir.

\begin{Lem} \label{diagrama}
Dado o diagrama
\begin{equation*}
 \xymatrix{
C(SS^2) \ar[d]_{p_2} \ar[r]^{p_1}    &          C(D \times S^1) \ar[d]^{\pi_1}      \\
C(D \times S^1) \ar[r]_{\pi_2}           &          C(S^1 \times S^1)                     } 
\end{equation*}
se as aplicações $p_1$ e $p_2$ forem restrições (inclusões) e
\begin{equation}\label{aplic}
 \begin{array}{lll}  \pi_1 :  &C(D \times S^1)  & \rightarrow C(S^1 \times S^1) \ \ \ \ \ \ \ \ \ \ \ \ \ \ \ \ \ \ \ \  e   \\
                                   & \ \ \ \ \ \ \ \varphi & \mapsto \pi_1 \circ \varphi = \varphi_{ \mid_{S^1 \times S^1}}  \\  \end{array} 
\end{equation}

\begin{equation} \label{aplic2}
 \begin{array}{lll}  \pi_2 :  &C(D \times S^1)  & \rightarrow C(S^1 \times S^1) \\
                                   & \ \ \ \ \ \ \ \varphi & \mapsto \pi_2 \circ \varphi (z,w) = \varphi_{ \mid_{S^1 \times S^1}} (z, -z^2 \bar{w}) ,\\  \end{array} 
\end{equation}
então 
$$C(SS^2) \simeq \{ (a_1,a_2) \in C(D \times S^1) \oplus C(D \times S^1) : a_1{\mid_{S^1 \times S^1}} = a_2{\mid_{S^1 \times S^1}}  \}.$$
\end{Lem}
\Dem
Vamos entender o porquê destas aplicações, na verdade o porquê de $\pi_2$.

Tomemos as cartas dadas pelas projeções estereográficas $\chi_N : S^2 - \{ PN \} \rightarrow \mathbb{R}^2$ e $\chi_S : S^2 - \{ PS \} \rightarrow \mathbb{R}^2$, mais precisamente,
 
$$\begin{array}{lll}  \chi_S :  & S^2 - \{ (0,0,-1) \}  & \rightarrow \mathbb{R}^2 \ \ \ \ \ \ \ \ \ \ \ \ \ \ \ \ \ \ \ \ \ \ \ \ \ \ \ \ \ \ \ \ \ \ \ \ \ \ \ \textrm{e}   \\
                                  & (x_1, x_2, x_3) & \mapsto \left( \frac{x_1}{(1+x_3)} , \frac{x_2}{(1+x_3)} \right) = (S_1 , S_2)\\  \end{array} $$  

$$\begin{array}{lll}  \chi_N :  & S^2 - \{ (0,0,1) \}  & \rightarrow \mathbb{R}^2   \\
                               & (x_1, x_2, x_3) & \mapsto \left( \frac{x_1}{(1-x_3)} , \frac{x_2}{(1-x_3)} \right) = (N_1, N_2) \\  \end{array}. $$ 

Consequentemente, no fibrado tangente, temos

$$\begin{array}{lll}  T \chi_S :  & T(S^2 - \{ (0,0,-1) \})  & \rightarrow \mathbb{R}^2 \times \mathbb{R}^2 \ \ \ \ \ \ \ \ \ \ \ \ \ \ \ \ \ \ \ \  \textrm{e}   \\
                                  &  T_x(S^2 - \{ (0,0,-1) \}) \ni \lambda & \mapsto \left( \chi_S(x) , \xi_1 , \xi_2 \right) \\  \end{array} $$  

$$\begin{array}{lll}  T \chi_N :  & T( S^2 - \{ (0,0,1) \})  & \rightarrow \mathbb{R}^2 \times \mathbb{R}^2  \\
                               &  T_x(S^2 - \{ (0,0,1) \}) \ni \lambda & \mapsto \left( \chi_N(x) , \eta_1 , \eta_2 \right) \\  \end{array} . $$ 

Para entendermos o comportamento da aplicação 
$$(\chi_N, \eta_1, \eta_2) \mapsto (\chi_S, \xi_1, \xi_2),$$ basta observar que a restrição $K : (\chi_N) = (N_1, N_2) \mapsto (\chi_S) = (S_1 , S_2)$ é dada por
$$\begin{array}{llll}  K : & \mathbb{R}^2 - \{ (0,0)\}  & \rightarrow  S^2 - \{ (0,0,1) , (0,0, -1)\} & \rightarrow  \mathbb{R}^2 - \{ (0,0)\}   \\
                             &  (N_1, N_2)   &  \mapsto \left( \frac{2N_1}{1+N_1^2+N_2^2}, \frac{2N_2}{1+N_1^2+N_2^2}, \frac{N_1^2 + N_2^2 - 1}{1+N_1^2+N_2^2} \right) & \mapsto \left( \frac{N_1}{N_1^2 + N_2^2} , \frac{N_2}{N_1^2 + N_2^2} \right) \\  \end{array} .$$
Portanto 
$\frac{\partial (K_1, K_2)}{\partial (N_1, N_2)} = \frac{1}{(N_1^2 + N_2^2)^2} \left( \begin{array}{ll}  - N_1^2 + N_2^2 & -2 N_1 N_2 \\ -2 N_1 N_2 &  N_1^2 - N_2^2 \\  \end{array} \right).$

Assim , para $(z, w) \in \mathbb{C}^2$, com $z = (N_1, N_2)$ e $w = (\eta_1, \eta_2)$, temos
$$T(z,w) =  T\chi_S \circ (T\chi_N)^{-1}(N_1, N_2, \eta_1, \eta_2) = \left(\frac{N_1}{N_1^2 + N_2^2} , \frac{N_2}{N_1^2 + N_2^2}, \frac{\partial (K_1, K_2)}{\partial (N_1, N_2)} \left( \begin{array}{ll} \eta_1 \\ \eta_2 \end{array} \right) \right),$$
que resumidamente é $T(z,w) = \left( \frac{z}{|z|^2} , \frac{-z^2 \bar{w}}{|z|^4} \right)$, mas quando restrito a $S^1 \times S^1$ temos $T(z,w) = (z, -z^2 \bar{w})$.

Isto explica o fato de $\pi_2$ ser dada por
$$ \begin{array}{lll}  \pi_2 :  &C(D \times S^1)  & \rightarrow C(S^1 \times S^1) \\
                                   & \ \ \ \ \ \ \ \varphi & \mapsto f_2 \circ \varphi (z,w) = \varphi_{ \mid_{S^1 \times S^1}} (z, -z^2 \bar{w}) .\\  \end{array} $$
\fimdem

Utilizando agora o Teorema \ref{MVH} e as informações do Lema \ref{diagrama} obtemos a seqüência exata cíclica

\begin{equation} \label{CSS2}
 \xymatrix{
K_0(C(SS^2))  \ar[rr]^-{(p_{1_*},p_{2_*})}  & \ \ &  K_0(C(D \times S^1))^2 \ \ \ar[rr]^-{\pi_{2_*} - \pi_{1_*}}     & \ \ &   K_0(C(S^1 \times S^1))    X
 \ar[dd]       \\
& \\
K_1(C(S^1 \times S^1)) \ar[uu]^{}  & \ \ &  K_1(C(D \times S^1))^2 \ar[ll]_-{\pi_{2_*} - \pi_{1_*}}  & \ \  & K_1(C(SS^2) \ar[ll]_-{(p_{1_*},p_{2_*})}                     } 
\end{equation}


Como podemos perceber, para fazer uso desta seqüência precisaremos antes saber quem são $K_i(C(D \times S^1))$ e $K_i (C(S^1 \times S^1))$, para $i=0$ e $i=1$, e seus geradores. E para realizarmos estes cálculos utilizaremos alguns resultados de K-teoria apresentados no segundo apêndice deste trabalho. 

Lembrando que $C(D \times S^1) \simeq C(S^1, C(D))$ e $C(S^1 \times S^1) \simeq C(S^1, C(S^1))$, podemos então construir as seqüências exatas curtas
$$ 0 \longrightarrow SC(D) \longrightarrow C(S^1, C(D)) \longrightarrow C(D) \longrightarrow 0$$
e
$$ 0 \longrightarrow SC(S^1) \longrightarrow C(S^1, C(S^1)) \longrightarrow C(S^1) \longrightarrow 0$$
e por \ref{cs1a} temos
$$K_i(C(D \times S^1)) = K_{1-i}(C(D)) \oplus K_i(C(D)) \ \ \ \textrm{ e }$$ 
$$K_i(C(S^1 \times S^1)) = K_{1-i}(C(S^1)) \oplus K_i(C(S^1)).$$

Dos exemplos \ref{cded} e \ref{Ks1}  e dos isomorfismos acima temos: 

\begin{itemize}
\item $K_0(C(D \times S^1)) = [\mathfrak{1}]_0 \mathbb{Z}$ 
\item $K_1(C(D \times S^1)) = [\mathfrak{w}]_1 \mathbb{Z}$
\item $K_0(C(S^1 \times S^1)) = [\mathfrak{1}]_0 \mathbb{Z} \oplus [\theta_{C(S^1)}(\mathfrak{w})]_0 \mathbb{Z}$ 
\item $K_1(C(S^1 \times S^1)) = [\mathfrak{z}]_1 \mathbb{Z}\oplus [\mathfrak{w}]_1 \mathbb{Z}$.
\end{itemize}

Com isto o diagrama (\ref{CSS2}) pode ser visto como

$$ \xymatrix{
K_0(C(SS^2))  \ar[rr]^-{(p_{1_*},p_{2_*})_0}  & \ \ &  \mathbb{Z} [\mathfrak{1}] \oplus \mathbb{Z} [\mathfrak{1}] \ \ \ar[rr]^-{(\pi_{2_*} - \pi_{1_*})_0}     & \ \ &   \mathbb{Z} [\mathfrak{1}] \oplus \mathbb{Z} [\theta_{C(S^1)}(\mathfrak{w})]  \ar[dd]     \\
& \\
\mathbb{Z} [\mathfrak{z}] \oplus \mathbb{Z} [\mathfrak{w}]  \ar[uu]^{}  & \ \ &  \mathbb{Z} [\mathfrak{w}] \oplus \mathbb{Z} [\mathfrak{w}]  \ar[ll]_-{(\pi_{2_*} - \pi_{1_*})_1}  & \ \  & K_1(C(SS^2). \ar[ll]_-{(p_{1_*},p_{2_*})_1}                     } $$
%
%

Pela definição dada nas equações (\ref{aplic}) e (\ref{aplic2}) do Lema \ref{diagrama} obtemos $$\pi_{1_*} (x,y) = (x,0) \textrm{   e   }   \pi_{2_*}(x,y) = (y,0),$$ ou seja, na linha superior do diagrama anterior temos 
\begin{equation} \label{pi0}
(\pi_{2_*} - \pi_{1_*})_0 (x,y) = (y - x,0). \end{equation}

Além disso, temos $\pi_{1_*} ([\mathfrak{w}]) = ([\mathfrak{w}]) $ e $\pi_{2_*}([\mathfrak{w}]) = [-\mathfrak{z}^2\bar{\mathfrak{w}}]$, ou seja, os geradores da ima\-gem serão $(0,1)$ e $(2,-1)$ e portanto para a linha inferior do diagrama anterior teremos o homomorfismo 
\begin{equation} \label{pi1} (\pi_{2_*} - \pi_{1_*})_1 (x,y) = (2y, -y-x). \end{equation}

Portanto veremos o diagrama anterior como 
$$ \xymatrix{
K_0(C(SS^2))  \ar[rr]^-{(p_{1_*},p_{2_*})_0}  & \ \ &  \mathbb{Z} [\mathfrak{1}] \oplus \mathbb{Z} [\mathfrak{1}] \ \ \ar[rr]^-{(x,y) \mapsto (y-x,0)}     & \ \ &   \mathbb{Z} [\mathfrak{1}] \oplus \mathbb{Z} [\theta_{C(S^1)}(\mathfrak{w})]  \ar[dd]^{\delta_0}     \\
& \\
\mathbb{Z} [\mathfrak{z}] \oplus \mathbb{Z} [\mathfrak{w}]  \ar[uu]^{\delta_1}  & \ \ &  \mathbb{Z} [\mathfrak{w}] \oplus \mathbb{Z} [\mathfrak{w}]  \ar[ll]_-{(x,y) \mapsto (2y, -y-x)}  & \ \  & K_1(C(SS^2). \ar[ll]_-{(p_{1_*},p_{2_*})_1}                     } $$
%

Note que o homomorfismo $(\pi_{2_*} - \pi_{1_*})_1$, isto é,  a aplicação $(x,y) \mapsto (2y, -y-x)$ é injetora, logo a imagem de $(p_{1_*},p_{2_*})_1$ é nula. E o núcleo de $(\pi_{2_*} - \pi_{1_*})_0$, dada por $(x,y) \mapsto ( y-x,0)$, é $\mathbb{Z}[(1,1)]$. Assim do diagrama anterior obtemos a seqüência exata
$$ 0 \longrightarrow \mathbb{Z}[(1,1)] \longrightarrow  K_1( C(SS^2))  \longrightarrow 0$$
e portanto o isomorfismo
$$K_1( C(S^*S^2)) \simeq \mathbb{Z}[(1,1)].$$

Temos ainda que o núcleo de $(\pi_{2_*} - \pi_{1_*})_0$ é isomorfo a $\mathbb{Z}[(1,1)]$, ou seja, a imagem de $(p_{1_*},p_{2_*})_0$ é $\mathbb{Z}[(1,1)]$. E a partir da imagem de $(\pi_{2_*} - \pi_{1_*})_1$ que é $2\mathbb{Z} \oplus \mathbb{Z}$, podemos obter o núcleo da aplicação $\gamma: \mathbb{Z}[\mathfrak{z}] \oplus \mathbb{Z}[\mathfrak{w}] \rightarrow K_0(C(SS^2))$  que é um grupo isomorfo a $(\mathbb{Z} \oplus \mathbb{Z} )/ (2\mathbb{Z} \oplus \mathbb{Z}) \simeq \mathbb{Z}_2[(1,0)]$. Com isto construímos a seqüência exata
$$ 0 \longrightarrow \mathbb{Z}_2 [(1,0)]  \longrightarrow  K_0( C(SS^2))  \longrightarrow  \mathbb{Z}[(1,1)] \longrightarrow 0$$
o que consequentemente nos leva a  
$$K_0( C(S^*S^2)) \simeq \mathbb{Z}[(1,1)] \oplus \mathbb{Z}_2[(1,0)].$$

\section{O caráter de Chern - Connes}

Como dito na introdução deste capítulo precisamos conhecer a K-Teoria da C$^*$-álgebra $\mathcal{A}$ gerada pelos operadores pseudodiferencias clássicos de ordem zero da esfera e para tal vamos utilizar a seqüência exata (\ref{seqS2}) e os grupos $K_0 (\mathcal{A}/\mathcal{K})$ e $K_1(\mathcal{A}/\mathcal{K})$ encontrados na seção anterior.

Novamente através da seqüência exata
\begin{equation*} 
 0 \longrightarrow \mathcal{K} \overset{\psi}{\longrightarrow} \mathcal{A} \overset{\pi}{\longrightarrow} \mathcal{A}/ \mathcal{K}  \longrightarrow 0 \ \ ,
\end{equation*}
como em \ref{seis}, a partir da funtoriedade de $K_0$ e $K_1$ e da periodicidade de Bott \ref{Bott}, obtemos a seqüência exata cíclica
$$ \xymatrix{
K_0(\mathcal{K})  \ar[rr]^-{K_0(\psi)}  & \ \ &  K_0(\mathcal{A}) \ \ \ar[rr]^-{K_0(\pi)}     & \ \ &   K_0( \mathcal{A}/\mathcal{K}) \ar[dd]^-{}     \\
& \\
K_1(\mathcal{A}/\mathcal{K})  \ar[uu]^-{\delta}  & \ \ &  K_1(\mathcal{A})  \ar[ll]_-{K_1(\pi)}  & \ \  & K_1(\mathcal{K}). \ar[ll]_-{K_1(\psi)}       } $$

Pelos cálculos da seção anterior temos 
$$K_0(\mathcal{A}/\mathcal{K}) \simeq K_0( C(SS^2)) \simeq \mathbb{Z}[(1,1)] \oplus \mathbb{Z}_2[(1,0)]$$
e
$$K_1(\mathcal{A}/\mathcal{K}) \simeq K_1( C(S^*S^2)) \simeq \mathbb{Z}[(0,1)].$$ 

Como em \cite{LM} a fórmula de Fedosov nos mostra que em toda variedade compacta, como é o nosso caso, existe um operador pseudodiferencial clássico de ordem zero com índice de Fredholm igual a $1$ (um), ou seja, a aplicação do índice $\delta$ é sobrejetora, logo
$$ \xymatrix{
\mathbb{Z}  \ar[rr]^-{K_0(\psi) = 0}  & \ \ &  K_0(\mathcal{A}) \ \ \ar[rr]^-{K_0(\psi)}     & \ \ &  \mathbb{Z}[(1,1)] \oplus \mathbb{Z}_2[(1,0)]  \ar[dd]^-{}     \\
& \\
\mathbb{Z}[(0,1)]  \ar[uu]^-{sobrej.}  & \ \ &  K_1(\mathcal{A})  \ar[ll]_-{K_1(\psi)}  & \ \  & \ 0. \ar[ll]_-{K_1(\pi)}       } $$

Assim, como já feito na seção anterior, podemos decompor este diagrama em duas seqüências exatas, e são elas
$$  0 \longrightarrow K_0(\mathcal{A}) \longrightarrow \mathbb{Z}[(1,1)] \oplus \mathbb{Z}_2[(1,0)]  \longrightarrow 0$$
isto porque o núcleo de $K_0(\psi) = 0$ e 
$$ 0 \longrightarrow K_1(\mathcal{A})  \longrightarrow 0$$
pois o homorfismo de grupos $\delta: \mathbb{Z} \rightarrow \mathbb{Z}$ é sobrejetor e conseqüentemente injetor, logo a imagem de $K_1(\psi) = 0$.

Destas duas seqüências temos
$$ K_0(\mathcal{A}) = \mathbb{Z}[(1,1)] \oplus \mathbb{Z}_2[(1,0)] \ \ \ \ \ \textrm{ e } \ \ \ \ \  K_1(\mathcal{A}) = 0.$$
Observe então que o caráter de Chern-Connes neste caso se resume ao caso par, ou seja, basta encontrarmos $Ch_\tau : K_0(\mathcal{A}) \rightarrow H^{\textrm{par}}_\mathbb{R} (SO(3))$, já que $K_1(\mathcal{A}) = 0$.

Antes de calcularmos o caráter de Chern-Connes como na Definição \ref{Chernpar}, cabe observar que qualquer homomorfismo $h:\mathbb{Z}_2 \rightarrow \mathbb{R}$ é nulo, já que 
$$0 = h(2) = h(1+1) = 2h(1).$$

Com isto o caráter de Chern-Connes para o C$^*$-sistema dinâmico em questão, isto é,  para $(\mathcal{A}=\overline{\Psi^0(S^2)}, SO(3), \alpha)$ é dado por
$$ ch_\tau : K_0(\mathcal{A}) \oplus K_1(\mathcal{A}) = \mathbb{Z}[(1,1)] \oplus \mathbb{Z}_2[(1,0)] \oplus 0 \rightarrow H^*_\mathbb{R}(SO(3))$$
e este não se anula apenas no gerador $[(1,1)]$ de $K_0(\mathcal{A})$ e portanto se resume a
$$ Ch_\tau : K_0(\mathcal{A}) = \mathbb{Z}[(1,1)]  \rightarrow \mathbb{R},$$
com $Ch_\tau ([1,1]_0) = \tau(\mathfrak{z})$ que é igual a medida de $S^2$.

\ \ \

\newpage

\ \ \ 

\newpage

\chapter{A K-teoria do fibrado das coesferas do Toro}
Dado o C$^*$-sistema dinâmico $(\mathcal{A}, S^1, \alpha)$, onde $\mathcal{A}$ é $\overline{\Psi_{cl}^0(T)}$ a C$^*$-álgebra gerada pelos operadores pseudodiferenciais clássicos de ordem zero do Toro e a ação $\alpha$ é a de conjugação pela translação em $S^1$, podemos repetir o raciocínio utilizado nos dois capítulos anteriores para calcular o caráter de Chern-Connes neste caso.

Note que assim como nos exemplos anteriores temos $\alpha$ contínua e $\Psi_{cl}^0(T)$ invariante pelo cálculo funcional holomorfo o que nos dá razão em dizer que estamos calculando o caráter de Chern-Connes para os operadores pseudodiferenciais clássicos de ordem zero do Toro.

E para calcularmos $K_0(\mathcal{A})$ e $K_1(\mathcal{A})$ neste caso, precisaríamos antes calcular a K-teoria do fibrado das coesferas do Toro e é por este motivo que apresentaremos neste capítulo o cálculo dos grupos de K-teoria do fibrado das coesferas do Toro e seus geradores.

De fato, para obter $K_0(\mathcal{A})$ e $K_1(\mathcal{A})$ seria necessário utilizar um raciocínio aná\-logo ao feito em (\ref{seqS1}) e (\ref{seqS2}) no qual $\mathcal{K} \subset \mathcal{A}$ e a sequência 
\begin{equation} \label{seqT}
 0 \longrightarrow \mathcal{K} \longrightarrow \mathcal{A} \overset{\pi}{\longrightarrow} \mathcal{A}/\mathcal{K}  \longrightarrow 0
\end{equation}
é exata.

Novamente pelo Teorema $2$ de \cite{TO} e seu corolário, existe um isomorfismo $\varphi$ entre o quociente $\mathcal{A}/\mathcal{K}$ e $C(S^*T) \overset{\psi}{\simeq} C(S^1 \times T)$. 

Utilizando a seqüência (\ref{seqT}) e o Teorema \ref{seis} temos a seqüência exata cíclica de seis termos
$$ \xymatrix{
K_0(\mathcal{K})  \ar[rr]^-{K_0(\psi)}  & \ \ &  K_0(\mathcal{A}) \ \ \ar[rr]^-{K_0(\pi)}     & \ \ &   K_0( \mathcal{A}/\mathcal{K} \simeq C(S^*T)) \ar[dd]^-{}     \\
& \\
K_1(\mathcal{A}/\mathcal{K} \simeq C(S^*T))  \ar[uu]^-{\delta}  & \ \ &  K_1(\mathcal{A})  \ar[ll]_-{K_1(\pi)}  & \ \  & K_1(\mathcal{K}). \ar[ll]_-{K_1(\psi)}       } $$

Assim como visto anteriormente, como $T$ é compacto, a aplicação $\delta$ é sobrejetora e portanto $K_1(\psi) = 0$, além disso $K_1(\mathcal{K}) = 0$ o que nos permite decompor o diagrama acima em duas seqüências exatas, são elas
$$0 \rightarrow K_0(\mathcal{A}) \rightarrow K_0(C(S^*T))\rightarrow 0 $$
e
$$0\rightarrow K_1(\mathcal{A}) \rightarrow K_1(C(S^*T)) \rightarrow \mathbb{Z} \rightarrow 0.$$
 
Através das duas seqüências exatas acima fica claro que um modo de encontrarmos $K_0(\mathcal{A})$ e $K_1(\mathcal{A})$ e seus geradores, é encontrando antes os grupos de K-teoria de $C(S^*T)$.

Vamos então calcular os grupos de K-teoria de $C(S^*T)$ e para isto pensemos em $C(T) = C(S^1 \times S^1) \simeq C(S^1, C(S^1))$, com $(z, w) \in S^1 \times S^1$, e tomemos a seguinte seqüência exata curta
$$ 0 \longrightarrow SC(S^1) \longrightarrow  C(T)  \longrightarrow  C(S^1) \longrightarrow 0.$$
Como a seqüência exata acima cinde, temos 
$$ K_0(C(T)) = K_0(C(S^1)) \oplus K_1(C(S^1)) = [1]_0 \mathbb{Z} \oplus [\theta_{C(S^1)} [\mathfrak{w}]_1]_0\mathbb{Z} $$
e
$$ K_1(C(T)) = K_1(C(S^1)) \oplus K_0(C(S^1)) = [\beta_{C(S^1)}[1]_0]_1 \mathbb{Z} \oplus [\mathfrak{w}]_1\mathbb{Z}, $$
lembrando que $\mathfrak{w}$ representa a função $\mathfrak{w}(z, w) = w$.

Pensando agora em $S^*T \simeq S^1 \times T = S^1 \times S^1 \times S^1$, $(\xi, z, w) \in S^*T$ e procedendo-se da mesma maneira, temos a seqüência exata cindida
$$ 0 \longrightarrow SC(T) \longrightarrow  C(S^1, C(T) )  \longrightarrow  C(T) \longrightarrow 0.$$
Logo 
\begin{equation*}
\begin{split}
 K_0(C(S^*T)) &= K_0(C(T)) \oplus K_1(C(T)) \\
 &= [1]_0 \mathbb{Z} \oplus [\theta_{C(S^1)} [\mathfrak{w}]_1]_0\mathbb{Z} \oplus [\theta_{C(T)}[\beta_{C(S^1)}[1]_0]_1]_0 \mathbb{Z} \oplus [\theta_{C(T)}[\mathfrak{w}]_1]_0\mathbb{Z}
\end{split}
\end{equation*}
e  
\begin{equation*}
\begin{split}
 K_1(C(S^*T)) &= K_1(C(T)) \oplus K_0(C(T)) \\
 &= [\beta_{C(S^1)}[1]_0]_1 \mathbb{Z} \oplus [\mathfrak{w}]_1\mathbb{Z} \oplus  [\beta_{C(T)}[1]_0]_1 \mathbb{Z} \oplus [\beta_{C(T)}[\theta_{C(S^1)} [\mathfrak{w}]_1]_0]_1\mathbb{Z}.
\end{split}
\end{equation*}

\ \ \ 

\newpage

\ \ \ 

\newpage

\addcontentsline{toc}{chapter}{Apêndices}
\appendix
\chapter{ Representações em álgebras de Lie }
\vspace{0.75 cm}
Nossa intenção neste apêndice é a de mostrar alguns resultados básicos sobre derivação e representações em álgebras de Lie que são utilizados no primeiro capítulo  desta tese. Um exemplo disto é o de que dado um C$^*$-sistema dinâmico $(A,G, \alpha)$, a definição dada em Definição \ref{Der} faz sentido, isto é, o de que a aplicação $\delta$ é de fato uma representação de $\mathfrak{g}$, álgebra de Lie do grupo $G$, sobre a álgebra de Lie das derivações sobre $A^\infty$. 

Gostaríamos aqui de agradecer novamente o inestimável auxílio dado por Daniel V. Tausk que produziu grande parte do texto e das explicações aqui expostas, assim como sua colaboração em outros pontos deste trabalho.

Sejam $G$ um grupo de Lie, $A$ um espaço de Banach (que em nosso caso será uma C$^*$-álgebra) e $\alpha: G \rightarrow GL(A)$ um homomorfismo (lembrando que $GL(A)$ é o grupo dos isomorfismos lineares contínuos de A e que contém $Aut(A)$ o grupo dos $*$-automorfismos da C$^*$-álgebra $A$). 

O vetor $x \in A$ é de classe $C^k$ ($0 \leqslant k \leqslant \infty$) se a aplicação $\alpha(x): g \mapsto \alpha_g(x)$ é  de classe $C^k$. E assim, como na Definição \ref{Der}, se $x$ é no mínimo de classe $C^1$ e $X \in \mathfrak{g}$ definimos 
$$\delta_X(x) = \frac{d}{dt}[\alpha_{exp(tX)}(x)]|_{t=0}.$$

A proposição a seguir é basicamente a demonstração do teorema de G\aa rding, isto é, a de que a *-álgebra de Fréchet $A^\infty$ é densa (na norma) em $A$, para maiores detalhes vide Teorema $1$ do capítulo $5$ de \cite{DA}.

\begin{Prop} \label{Garding}
Sejam $f \in C_c(G)$ e  $\nu$ uma medida de Haar (invariante à esquerda) sobre $G$. O operador $\alpha_f: A \rightarrow A$ dado por $\alpha_f(a) = \int_G f(x)\alpha_x(a) d \nu (x)$ satisfaz as seguintes afirmações.
\begin{enumerate}
\item Se $f \in C_c^\infty(G)$, então $\alpha_f(a) \in A^\infty$, $\forall a \in A$.
\item Dada uma seqüência $(f_i)_{i \in \mathbb{N}} \subset C_c^\infty(G)$ tal que: \\
($i$) $\int_G f_i(x) d \nu(x) = 1$, $\forall i \in \mathbb{N}$; \\
($ii$) $f_i \geqslant 0$, $\forall i \in \mathbb{N}$  e  \\
($iii$) para qualquer vizinhança $U$ da unidade de $G$, temos $supp \ f_i \subset U$ exceto para um número finito de índices. \\
Então $\alpha_{f_i}(a)$ converge para $a$,  para qualquer $a \in A$.
\end{enumerate}
\end{Prop}
\Dem (1)
Seja $g$ um elemento de $G$, como a medida de Haar é invariante à esquerda obtemos 
$$\alpha_g \alpha_f(a) = \int_G f(x) \alpha_{gx} (a) d \nu (x) = \int_G f(g^{-1}x) \alpha_{x} (a) d \nu (x).$$ 

Assim para $X \in \mathfrak{g}$ e $g_t = exp(tX)$, com $t \in \mathbb{R}$, temos 
$$f(g_t^{-1}x) - f(x) = \int_0^1 \frac{d}{ds} f(g_{st}^{-1}x) ds = -t\int_0^1 Xf(g_{st}^{-1}x)ds , $$ vendo aqui $X$ como o campo de vetores invariante à direita sobre $G$ e que corresponde a transformação $x \mapsto g_tx$.

Portanto
\begin{equation*}
\begin{split}
\frac{\alpha_{g_t} - I}{t}\alpha_f(a) &= - \int_G \int_0^1 Xf(g^{-1}_{st}x)\alpha_x (a) ds d\nu(x) \\
&= - \int_G \int_0^1 Xf(x) \alpha_{g_{st}} \alpha_x (a) ds d\nu(x) \\
&= - \int_0^1  \alpha_{g_{st}} \int_G Xf(x) \alpha_x (a)  d\nu(x) ds \\
&= - \int_0^1 \alpha_{g_{st}}\alpha_{Xf}(a) ds 
\end{split}
\end{equation*}
e
$$\underset{t\rightarrow 0}{\lim} \frac{\alpha_{g_t} - I}{t}\alpha_f(a) = - \alpha_{Xf}(a) .$$

Provamos com isto que a aplicação $\alpha(\alpha_f(a)) : G \rightarrow A$, dada por $g \mapsto \alpha_g \alpha_f(a)$ é de classe $C^1$ e repetindo-se este mesmo argumento podemos provar que esta é uma aplicação de classe $C^\infty$.

(2) Dados $\epsilon >0$ e $a \in A$, podemos escolher uma vizinhança $U$ da unidade de $G$ tal que $||\alpha_g (a) - a|| \leqslant \epsilon$, para $g \in U$. Como $\alpha_{f_i}(a) - a = \int_G f_i(g)(\alpha_gf-f) d \nu(g)$, então se $supp \ f_i \subset U$, temos $||\alpha_{f_i}(a)-a|| \leqslant \epsilon$.
\fimdem

\begin{Def}
Dadas $\mathfrak{g}_1$ e $\mathfrak{g}_2$ álgebras de Lie dizemos que a aplicação $\varphi: \mathfrak{g}_1 \rightarrow \mathfrak{g}_2$ é um homomorfismo de álgebras de Lie se $\varphi$ é linear e preserva os colchetes de Lie, isto é, $\varphi([X,Y]) = [\varphi(X), \varphi(Y)]$, $\forall X,Y \in \mathfrak{g}_1$. 

Além disso, se $\mathfrak{g}_2$ for $GL(n, \mathbb{C})$, $GL(n, \mathbb{R})$, ou ainda $End(V)$, $V$ espaço vetorial, o homomorfismo $\varphi$ é chamado de representação.
\end {Def}

Dadas duas variedades diferenciáveis $M$ e $N$, $A$ um espaço de Banach e uma aplicação $f: M \times N \rightarrow A$ de classe $C^2$, denotaremos, para $(x_0, y_0) \in M \times N$, a diferencial da aplicação $x \mapsto f(x, y_0)$ em $x_0$ por
$$\partial_1 f(x_0,y_0) : T_{x_0}M \rightarrow A $$
e a diferencial da aplicação $y \mapsto f (x_0, y)$ no ponto $y_0$ por 
$$\partial_2 f(x_0,y_0) : T_{y_0}N \rightarrow A .$$

Apenas utilizando o Teorema de Schwarz em coordenadas locais em torno de $(x_0, y_0)$ podemos mostrar que para quaisquer $v \in T_{x_0}M$ e $w \in T_{y_0}N$ a diferencial da aplicação
$$ x \mapsto \partial_2f(x,y_0).w \in A$$
no ponto $x_0$ avaliada em $v$ coincide com a diferencial da aplicação  
$$ y \mapsto \partial_1f(x_0,y).v \in A$$
no ponto $y_0$ avaliada em $w$.

\begin{Lem}
Se $x \in A$ é de classe $C^k$, então:
\begin{enumerate}
\item $\forall g \in G$ o vetor $\alpha_g(x) \in A$ é de classe $C^k$;
\item se $k \geqslant 1$ e $X \in \mathfrak{g}$, então a aplicação $\lambda_X : g \mapsto \delta_X(\alpha_g(x)) \in A$ é de classe $C^{k-1}$;
\item se $k \geqslant 1$, então $\forall Y \in \mathfrak{g}$ o vetor $\delta_Y (x) \in A$ é de classe $C^{k-1}$;
\item se $k \geqslant 2$ e $X,Y \in \mathfrak{g}$, então $d\lambda_X(1).Y = \delta_X \delta_Y(x)$;
\item se $k \geqslant 1$ e $X \in \mathfrak{g}$, então a aplicação $\tau_X: (g,h) \mapsto \alpha_g(\delta_X (\alpha_h(x))) \in A$ é de classe $C^{k-1}$;
\item se $k \geqslant 2$ e $X,Y \in \mathfrak{g}$, então 
$$\delta_1 \tau_X(1,1).Y = \delta_Y (\delta_X (x)) \ \ \textrm{ e } \ \ \delta_2 \tau_X(1,1).Y = \delta_X (\delta_Y (x)).$$
\end{enumerate}
\end{Lem}
\Dem
(1) Consideremos a aplicação
\begin{displaymath}
\begin{array}{lll}
f:& G \times G & \rightarrow  A \\
& (g,h) & \mapsto \alpha_{hg}(x) = \alpha_h \alpha_g(x)\\ 
\end{array},
\end{displaymath}
 de classe $C^k$. Como por hipótese $x$ é de classe $C^k$, então para todo $g \in G$ a aplicação $h \mapsto \alpha_h(\alpha_g(x))$ é de classe $C^k$, logo $\alpha_g(x)$ é de classe $C^k$.

(2) Para $k \geqslant 1$ temos
$\partial_2f(g,1).X = d (\alpha_1(\alpha_g (x))).X = \delta_X (\alpha_g(x)) = \lambda_X(g)$, consequentemente $\lambda_X$ é de classe $C^{k-1}$.

(3) Vale também para todo $h \in G$ e $Y \in \mathfrak{g}$ que 
$$\partial_1f(1,h).Y = d (\alpha_h(\alpha_1 (x))).Y = \alpha_h(\delta_Y (x)),$$ ou seja, $\delta_Y(x)$ é de classe $C^{k-1}$.

(4) Se $k \geqslant 2$, utilizando a afirmação anterior a este lema aplicada a $f$ no ponto $(1,1) \in G \times G$ e nos vetores $X,Y \in \mathfrak{g}$, temos
\begin{equation*}
\begin{split}
d\lambda_X(1).Y &= d(\partial_2f(1,1).X).Y = d(\partial_1f(1,1).Y).X \\
 &= d (\delta_Y(x)).X = d(\alpha_1(\delta_Y(x))).X = \delta_X (\delta_Y(x)).
\end{split}
\end{equation*}

(5) Consideremos agora a aplicação
\begin{displaymath}
\begin{array}{lll}
F:& G \times G \times G & \rightarrow  A \\
& (g_1,g_2,g_3) & \mapsto \alpha_{g_1g_3g_2}(x)\\ 
\end{array},
\end{displaymath}
novamente $F$ é de classe $C^k$ e para $y = \alpha_{g_2}(x)$ temos
$$\partial_3 F(g_1, g_2, 1) = d(\alpha_{g_1}(\alpha_1(y))) = \alpha_{g_1}( \delta _X(y)) = \tau_X (g_1, g_2),$$ ou seja, $\tau$ é de classe $C^{k-1}$.

(6) E se $k \geqslant 2$ a derivada $\partial_1 (\tau_X(1,1))$ é calculada observando-se que a aplicação $g \mapsto \tau_X(g,1)$ coincide com $\alpha_1(z)$, para $z=\delta_X(x)$. Assim $$\partial_1(\tau_X(1,1)).Y = d(\alpha_1(z)).Y = \delta_Y(\alpha_1(z)) = \delta_Y(\delta_X(x)).$$ Da mesma forma, $g \mapsto \tau_X(1,g)$ coincide com $\lambda_X$ e 
$$\partial_2(\tau_X(1,1)).Y = d(\lambda_X(1)).Y = \delta_X(\delta_Y(x)).$$
\fimdem

\begin{Teo} \label{representacao}
Se $x \in A$ é de classe $C^2$, então $\forall X,Y \in \mathfrak{g}$ temos $$ \delta_{[X,Y]}(x) = \delta_X \delta_Y(x) - \delta_Y \delta_X(x) .$$
\end{Teo}
\Dem
Dados $g,h \in G$ temos $\alpha_g \alpha_h \alpha_{g^{-1}} (x) = \alpha_{i_g(h)}(x)$, onde $i_g$ é o automorfismo interno de $G$, isto é, $i_g(h) = ghg^{-1}$. Da mesma forma, se $y= \alpha_{g^{-1}}(x)$, temos $\alpha_g \circ \alpha(y) = \alpha(x) \circ i_g $ \footnote{a aplicação $\alpha(x): G \rightarrow A$ é aquela dada por $\alpha(x)(g) = \alpha_g(x)$.}. Pelo item (1) do lema anterior, $y$ é de classe $C^2$, desse modo diferenciando a igualdade anterior em $1 \in G$ de ambos os lados e avaliando em $Y \in \mathfrak{g}$ temos:
$$d(\alpha_g \alpha(y)(1)).Y = d(\alpha(x)i_g(1)).Y$$ 
$$ \alpha_g(\delta_Y(y))(x) = \delta_{Ad_g(Y)}(x),$$
para todo $g \in G$, onde $Ad_g = d (i_g(1)): \mathfrak{g} \rightarrow \mathfrak{g}$.

Diferenciando-se a igualdade acima em $g = 1$ e avaliando-a em $X \in \mathfrak{g}$, lembrando que $ \alpha_g (\delta_Y (\alpha_{g^{-1}} (x))) = \tau_Y(g, g^{-1})$, temos 
$$ d(\alpha_g (\delta_Y (\alpha_{g^{-1}} (x)))).X = d(\tau_Y(g, g^{-1})).X $$
 e em $g=1$ temos $d(\tau_Y(1, 1)).X $ que pelo lema anterior é de classe $C^1$ e a diferencial avaliada em $X$ é dada por 
$\partial_1 \tau_Y (1,1).X - \partial_2 \tau_Y(1,1).X = \delta_X \delta_Y (x)-\delta_Y \delta_X (x)$.

Agora diferenciando-se o lado direito da igualdade temos $Z \mapsto \delta_Z(x) \in A$ linear e contínua, pois coincide com $d(\alpha(x)(1))$. Além disso, a diferencial da aplicação $g \mapsto Ad_g(Y) \in \mathfrak{g}$ no ponto $g=1$ e avaliada em $X$ é igual a $[X,Y]$, ou seja, 
$$d(\delta_{Ad_1(Y)}(x)).X = \delta_{[X,Y]}(x).$$
\fimdem

Com isto provamos que a aplicação $\delta$, como na Definição \ref{Der}, é de fato uma representação de $\mathfrak{g}$ na álgebra de Lie das derivações de $A^\infty$. 

\ \ \

\newpage

\chapter{ Resultados sobre K-teoria de C$^*$-álgebras }
Este apêndice dedica-se a uma breve introdução e à apresentação de alguns resultados básicos, utilizados neste trabalho, sobre a K-teoria de C$^*$-álgebras. Os resultados aqui apresentados são clássicos e podem ser vistos com mais detalhes em livros como \cite{RLL}, \cite{WO} e \cite{BB}, que nesta seqüência, aumentam o grau de complexidade e diminuem o de didática.

\begin{Def}
Uma C$^*$-álgebra $A$ é uma $*$-álgebra de Banach onde $\| a^*a \| = \|a\|^2$, para todo $a \in A$. Se $A$ possui unidade, $1_A$, então dizemos que $A$ é uma C$^*$-álgebra com unidade. 
\end{Def}

Alguns exemplos bastante conhecidos e utilizados de C$^*$-álgebras são

\begin{Ex} \label{BH}
Seja $\mathcal{H}$ um espaço de Hilbert de dimensão finita ou infinita, então $\mathcal{L}(\mathcal{H})$, conjunto dos operadores limitados de $\mathcal{H}$, é uma C$^*$-álgebra quando conside\-radas as operações
\begin{equation*}
\begin{split}
 (T+S)(h) &= T(h) + S(h) \\
 (TS)(h) &= T\big( S(h) \big) \\
 (\lambda T)(h) &= \lambda T(h), \\
\end{split}
\end{equation*}
a norma usual de operadores
$$ \| T \| = \sup \{ \| Th \| : h \in H, \ \|h\| \leqslant 1 \},$$
e a involução obtida através do produto interno,
$$ \langle Th, k \rangle \ = \ \langle h, T^*k \rangle \ \ \  para \ h \ e \ k \ em \ H.$$
\end{Ex}

\begin{Ex}
O conjunto das matrizes quadradas $n \times n$ a valores complexos, $M_n(\mathbb{C})$, também é uma C$^*$-álgebra, basta identificarmos este espaço com $\mathcal{L}(\mathbb{C}^n)$ e temos um caso particular do exemplo anterior.  
\end{Ex}

\begin{Ex}
As subálgebras fechadas e auto-adjuntas de $\mathcal{L}(\mathcal{H})$, $\mathcal{H}$ um espaço de Hilbert, como no exemplo (\ref{BH}) também são C$^*$-álgebras.
\end{Ex}

\begin{Ex}
Dada uma C$^*$-álgebra A, o conjunto das matrizes $M_n(A)$, munido das operações usuais também é uma C$^*$-álgebra. 
Pelo Teorema de Gelfand-Naimark, vide Teorema 1.1.3 de \cite{RLL}, $A$  é isomorfa a uma C$^*$-subalgebra de $\mathcal{L}(\mathcal{H})$, para algum espaço de Hilbert $\mathcal{H}$. Assim $M_n(A)$ é uma C$^*$-subalgebra de $\mathcal{L}(\mathcal{H}^n)$.  
\end{Ex}

Dada uma $*$-álgebra $A$ qualquer, com ou sem unidade, podemos associar a esta uma $*$-álgebra $\widetilde{A}$ que possua uma unidade. Isso pode ser feito se construirmos $\widetilde{A}$, chamada de unitização de $A$, da seguinte forma: 
\begin{equation*}
\begin{split}
\widetilde{A} &= A \oplus \mathbb{C} \\
(a,z)(b,w) &= (a b + z b + w a, zw) \\
(a,z)^* &= (a^*,\overline{z})\\
\end{split}
\end{equation*}
e teremos como unidade de $\widetilde{A}$ o elemento $1_{\widetilde{A}} = (0,1)$. Além disso, teremos uma identificação, através de um homomorfismo injetivo, de $A$ com um ideal maximal de $\widetilde{A}$ através dos elementos $(a,0)$, para todo $a \in A$. 

Como $A$ pode ser identificado com uma subálgebra fechada de $\widetilde{A}$, então $\widetilde{A}$ será uma $*$-álgebra de Banach se $A$ de fato for uma $*$-álgebra de Banach. E por fim, se $A$ é uma C$^*$-álgebra podemos definir em $\widetilde{A}$ uma norma que a torna uma C$^*$-álgebra, vide Teorema 1.15 de \cite{DPD}.

Dadas as $*$-álgebras $A$ e $B$ podemos também estender os $*$-homomorfismos $\phi: A \rightarrow B$ para $*$-homomorfismos $\widetilde{\phi} : \widetilde{A} \rightarrow \widetilde{B}$, dados por 
$\widetilde{\phi} (x + z 1_{\widetilde{A}}) = \phi(x) + z 1_{\widetilde{B}}.$ 
 
Como pode ser visto em 1.1.6 de \cite{RLL} a associação entre $A$ e $\widetilde{A}$ nos permite construir a seqüência exata curta cindida, dada por
\begin{equation} \label{unit}
0 \rightarrow A \overset{i}{\hookrightarrow} \widetilde{A} \underset{\lambda}{\overset{\pi}{\rightleftharpoons }}  \mathbb{C} \rightarrow 0
\end{equation}
com $\pi$ a aplicação quociente e $\lambda(z) = z 1_{\widetilde{A}}$.
Além disso, se $A$ possui unidade a aplicação $(a,z) \mapsto a + z ( 1_{\widetilde{A}} - 1_A)$ é um isomorfismo  entre $ A \oplus \mathbb{C}$ e $\widetilde{A} $.

\begin{Def} \label{proj}  
Sendo $\mathcal{P}(A)$ o conjunto das projeções de uma C$^*$-álgebra $A$, isto é, o conjunto dos elementos de $p \in A$ tais que  $ p = p^* = p^2$, definimos 
$$ \mathcal{P}_n(A) = \mathcal{P} (M_n(A)) \ \ \ \textrm{ e } \ \ \ \mathcal{P}_\infty (A) = \overset{\infty}{\underset{n=1}{\cup }}\mathcal{P}_n(A) .$$

Dadas as projeções $p$ e $q$ em $\mathcal{P}_\infty (A)$, mais precisamente $p \in \mathcal{P}_m (A)$ e $q \in \mathcal{P}_n (A)$, definimos a relação de equivalência $p \sim_0 q$ se existe $u \in M_{mn}(A)$ tal que $p = uu^*$ e $q = u^*u$. 

Definimos também a operação 
$p \oplus q = \left( \begin{array}{ll}  p & 0 \\ 0 & q  \end{array} \right)$
em $\mathcal{P}_\infty (A)$.
\end{Def}

Utilizando a definição acima podemos construir o semigrupo abeliano e com unidade $(\mathcal{P}_0(A), +)$, onde 
$$ \mathcal{P}_0 (A) = \mathcal{P}_\infty (A)/ \sim_0 \ \ \ \textrm{ e } \ \ \ [p]_0 + [q]_0 =[ p \oplus q]_0.$$

Dado qualquer semigrupo abeliano $(S,\oplus)$ existe uma construção bastante co\-nhecida, chamada construção de Grothendieck, de um grupo abeliano a partir deste semigrupo. Definimos a relação que facilmente se prova ser de equivalência $\sim $ em $S \times S$, dada por $(x_1, y_1) \sim (x_2, y_2)$ se existe $z \in S$ tal que $x_1 \oplus y_2 \oplus z = x_2 \oplus y_1 \oplus z $. Dessa forma definimos o grupo de Grothendieck como sendo $G_S = S \times S / \sim$, cuja operação $+$ é dada por $[x_1,y_1] + [x_2, y_2] = [x_1 \oplus x_2 , y_1 \oplus y_2]$, onde $[x,y]$ designa a classe de equivalência de $(x,y) \in S \times S$. 

Não é difícil provar que $(G_S, +)$ construido acima é de fato um grupo e também que é abeliano.
Além disso, para $y \in S$ existe uma aplicação, conhecida como aplicação de Grothendieck, que independe do $y$ escolhido e que é dada por  
\begin{equation*}
\begin{split} 
\gamma_S : S & \rightarrow G_S \\
     x & \mapsto [x \oplus y, y] .\\ 
\end{split}
\end{equation*}

Estamos prontos portanto para definir o grupo de K-teoria $K_0(A)$ de uma C$^*$-álgebra $A$ \emph{com unidade}.

\begin{Def}
Dada uma C$^*$-álgebra $A$ com unidade, podemos definir o grupo de K-teoria $K_0(A)$ com sendo o grupo de Grothendieck obtido a partir do semigrupo abeliano e com unidade $(\mathcal{P}_0 (A), +)$, onde por abuso de notação $[p]_0 = \gamma([p]_0)$, para $\gamma: \mathcal{P}_0 \rightarrow K_0(A)$ a aplicação de Grothendieck.
\end{Def}

Existe uma representação padrão, muito utilizada, para os elementos de $K_0(A)$. Tal representação é apresentada na proposição a seguir e sua demonstração pode ser encontrada em $3.1.7$ de \cite{RLL}.

\begin{Prop}
Seja $A$ uma C$^*$-álgebra com unidade, então
\begin{equation*}
\begin{split}
K_0(A) &= \{ [p]_0 - [q]_0 : p, q \in \mathcal{P}_\infty (A) \} \\
 &= \{ [p]_0 - [q]_0 : p, q \in \mathcal{P}_n (A) \}. \\
\end{split}
\end{equation*}
\end{Prop}

A definição anterior nos permite entender $K_0$ como um funtor entre a categoria das C$^*$-álgebras com unidade e a categoria dos grupos abelianos, vide $3.2.4$ de \cite{RLL}. Além disso, poderíamos utilizar a definição acima também para C$^*$-álgebras \emph{sem unidade}, mas isto impediria o funtor $K_0$ de ser meio-exato
\footnote{Um funtor F é dito meio-exato se dada uma seqüência exata curta $ 0 \rightarrow  A \overset{\phi}{\rightarrow} B \overset{\psi}{\rightarrow} C \rightarrow 0$, temos a seqüência exata
$ F(A) \overset{F(\phi)}{\rightarrow} F(B) \overset{F(\psi)}{\rightarrow} F(C)$. }, 
como pode ser visto no exemplo $3.3.9$ também de \cite{RLL}.

Dado um $*$-homomorfismo de C$^*$-álgebras $\phi : A \rightarrow B$ utilizamos a funtoriedade de $K_0$ para definir o homomorfismo de grupos 
\begin{equation*}
\begin{split} 
K_0(\phi): K_0(A) & \rightarrow K_0(B) \\
     [p]_0 & \mapsto [\phi(p)]_0.\\ 
\end{split}
\end{equation*}

Já para C$^*$-álgebras $A$ sem unidade definimos o grupo $K_0(A)$, como sendo o núcleo da aplicação  $K_0(\pi) : K_0(\widetilde{A}) \rightarrow K_0(\mathbb{C})$, dada por $K_0(\pi)([p]_0) = [\pi(p)]_0$ que pode ser melhor entendida através da seqüência exata \ref{unit}, ou resumidamente 
$$ K_0(A) = N (K_0(\pi)).$$

E como pode ser visto em $4.1.3$, $4.3.2$ e $4.3.3$ de \cite{RLL} esta definição de $K_0$, para uma C$^*$-álgebra qualquer nos permite provar que além de funtor, $K_0$ também é meio-exato e exato cindido, isto é, leva seqüências exatas cindidas de C$^*$-álgebras em seqüências exatas cindidas de grupos abelianos.

Além disso, o funtor $K_0$ é contínuo com respeito a limite indutivo e é estável, isto é $K_0(\underset{\rightarrow }{\lim}A_n) = \underset{\rightarrow }{\lim} (K_0(A_n))$ e $K_0(A) \simeq K_0(A \otimes \mathcal{K})$ \footnote{Lembrando que $\mathcal{K}$ é a C$^*$-álgebra dos operadores compactos sobre um espaço de Hilbert separável e de dimensão infinita.}.

Assim como para C$^*$-álgebras com unidade, existe uma representação padrão para os elementos de $K_0(A)$, com $A$ uma C$^*$-álgebra qualquer.

\begin{Prop} \label{padrao}
Sejam $A$ uma C$^*$-álgebra qualquer e $s = \lambda \circ \pi : \widetilde{A} \rightarrow \widetilde{A}$, obtida através de \ref{unit}, então
\begin{equation*}
\begin{split}
K_0(A) &= \{ [p]_0 - [s(p)]_0 : p \in \mathcal{P}_\infty (\widetilde{A}) \} \\
\end{split}
\end{equation*}
\end{Prop}

Feito isto temos a definição, o retrato e algumas propriedades básicas do grupo $K_0(A)$ e do funtor $K_0$, para $A$ uma C$^*$-álgebra qualquer, mas antes de darmos alguns exemplos do cálculo de alguns desses grupos definiremos o outro grupo de K-teoria de uma C$^*$-álgebra $A$, em outros termos, $K_1(A)$.

\begin{Def} \label{unitario}
Sendo $\mathcal{U}(A)$ o grupo dos elementos unitários de uma C$^*$-álgebra $A$ com unidade, isto é, o conjunto dos elementos $u \in A$ tais que $u^*u = uu^* = 1_A$,  definimos 
$$\mathcal{U}_n(A) = \mathcal{U} (M_n(A)) \ \ \ \textrm{ e } \ \ \ \mathcal{U}_\infty (A) = \overset{\infty}{\underset{n=1}{\cup }}\mathcal{U}_n(A).$$

Definimos em $\mathcal{U}_\infty (A)$ a operação 
$u \oplus v = \left( \begin{array}{ll}  u & 0 \\ 0 & v  \end{array} \right).$

E dados os unitários $u$ e $v$ em $\mathcal{U}_\infty (A)$, ou mais precisamente, $u \in \mathcal{U}_n (A)$ e $v \in \mathcal{U}_m (A)$ definimos a relação de equivalência $ u \sim_1 v$ se existe $k \in \mathbb{N}$ tal que $u \oplus 1_{k-n} \sim_h v \oplus 1_{k-m}$ em $\mathcal{U}_k(A)$, onde $\sim_h$ é a usual homotopia por caminhos. 
\end{Def}

Estamos prontos portanto para definir o grupo $K_1(A)$ para uma C$^*$-álgebra $A$ qualquer.

\begin{Def} 
Dada uma C$^*$-álgebra $A$ qualquer definimos $$K_1(A) = \mathcal{U}_\infty (\widetilde{A}) / \sim_1.$$ Denotaremos por $[u]_1$ a classe de equivalência que contém o elemento $u \in \mathcal{U}_\infty (A)$ e a operação do grupo como sendo $[u]_1 + [v]_1 = [u \oplus v]_1$, para $u, v \in \mathcal{U}_\infty (A)$.
\end{Def}

Assim como $K_0$, $K_1$ também apresenta uma representação padrão que é dada pela proposição a seguir e que é uma consequência imediata da definição anterior.

\begin{Prop}
Dada uma C$^*$-álgebra $A$, então 
$$K_1(A) = \{ [u]_1 : u \in \mathcal{U}_\infty (\widetilde{A}) \} ,$$
com $[u]_1 + [v]_1 = [u \oplus v]_1$ e $[1]_1 = 0$.
\end{Prop}

Vale ressaltar, vide 8.1.6 de \cite{RLL}, que dada uma C$^*$-álgebra $A$ com unidade temos $K_1(A) \simeq \mathcal{U}_\infty (A) / \sim_1$.

Além disso, $K_1$ também é um funtor da categoria das C$^*$-álgebras na categoria dos grupos abelianos, vide $8.2$ de \cite{RLL}, também é meio exato, exato cindido, contínuo para limite indutivo e estável e as demonstrações destes fatos podem se encontradas em $8.2.4$, $8.2.5$, $8.2.7$ e $8.2.8$ de \cite{RLL}, respectivamente. 

A funtoriedade de $K_0$ e $K_1$, nos permite provar uma série de resultados comumente utilizados, como por exemplo.

\begin{Prop}
Dadas as C$^*$-álgebras $A$ e $B$, temos
$$ K_i (A \oplus B) \simeq K_i(A) \oplus K_i(B) \ \ \ \ \ \ \ \ \  i =0 \textrm{ ou } 1.$$
\end{Prop}

Estamos prontos agora para darmos alguns exemplos dos grupos $K_0(A)$ e $K_1(A)$, para algumas C$^*$-álgebras que serão utilizadas neste trabalho.

\begin{Ex} \label{KdeC}
Para a C$^*$-álgebra $A = M_n(\mathbb{C})$, particularmente $\mathbb{C} = M_1(\mathbb{C})$, temos 
$$K_0(A) = \mathbb{Z} \ \ \ \textrm{ e } \ \ \ K_1(A) = 0 .$$ 

Temos $K_1(A) = 0$ pois $A$ é conexo e para $Tr$ o traço usual das matrizes temos o isomorfismo $ K_0 (Tr)$ entre $K_0(A)$ e $\mathbb{Z}$, como pode ser visto em $8.1.8$ e $3.3.2$ de \cite{RLL} respectivamente.
\end{Ex}

Cabe ressaltar que dado um traço $\tau: A \rightarrow \mathbb{C}$ de uma C$^*$-álgebra $A$ qualquer, temos $\tau(p) = \tau(q)$, se $p \sim_0 q$. Assim $K_0(\tau)([p]_0) = \tau(p)$, para $p \in \mathcal{P}_\infty (A)$. 

Além disso, se $\tau$ é um traço positivo, isto é, $\tau(a) \geqslant 0$, para todo elemento positivo de $a \in A$, então $K_0(\tau)([p]_0) = \tau(p) \geqslant 0$, ou seja, $K_0(\tau) : K_0(A)\rightarrow \mathbb{R}$.

Da estabilidade dos funtores $K_0$ e $K_1$ e do exemplo anterior, temos 

\begin{Ex} \label{compacto}
Para $A = \mathcal{K}$, álgebra dos operadores compactos de um espaço de Hilbert separável, temos
$$K_0(A) = \mathbb{Z} \ \ \ \textrm{ e } \ \ \ K_1(A) = 0 .$$ 
\end{Ex}

\begin{Ex}
Seja $A = \mathcal{L}(\mathcal{H})$ a C$^*$-álgebra definida em \ref{BH}, para $\mathcal{H}$ separável ou não, temos 
$$ K_0(A) =  0 \ \ \ \textrm{ e } \ \ \  K_1(A) = 0.$$
\end{Ex}

\begin{Ex}
Seja $A = C(X)$, a C$^*$-álgebra das funções contínuas de um espaço contrátil $X$, temos
$$K_0(A) = \mathbb{Z} \ \ \ \textrm{ e } \ \ \ K_1(A) = 0.$$
\end{Ex}

A demonstração do exemplo anterior pode ser vista em $3.3.6$ e $8.2$ de \cite{RLL} e como caso particular deste temos alguns exemplos bastante utilizado neste trabalho, como 
\begin{Ex} \label{cded}
$$ K_0(C(D)) = K_0(C ( \{ pto \} ) ) = K_0(C([0,1])) = \mathbb{Z}$$
 e 
$$K_1(C(D)) = K_1(C ( \{ pto \} ) ) = K_1(C([0,1])) = 0$$
onde $D$ é o disco unitário, isto é, $D = \{ z \in \mathbb{C} : |z| \leqslant 1\}$ e o gerador de $\mathbb{Z}$ é $[\mathfrak{1}]_0$, onde $\mathfrak{1}$ é a unidade de $C(X)$. 
\end{Ex}

\begin{Def} \label{indice}
Dada uma seqüência exata curta de C$^*$-álgebras
$$ 0 \rightarrow A \overset{\phi}{\rightarrow} B \overset{\psi}{\rightarrow} C \rightarrow 0$$
existe $*$-homomorfismo natural $\delta_1 : K_1(C) \rightarrow K_0(A)$, chamado de aplicação do índice que torna a seqüência 
$$K_1(A)  \overset{K_1(\phi)}{\longrightarrow}     K_1(B) \overset{K_1(\psi)}{\longrightarrow} K_1(C) \overset{\delta_1}{\longrightarrow}  K_0(A)  \overset{K_0(\phi)}{\longrightarrow}  K_0(B)  \overset{K_0(\psi)}{\longrightarrow}  K_0(C) $$
exata. Além disso, existe uma representação padrão para tal aplicação dada por
\begin{equation*} 
\delta_1([u]_1) = [p]-[p_n]
\end{equation*}
onde $u \in \mathcal{U}_n(\widetilde{C})$, $v \in \mathcal{U}_{2n}(\widetilde{B})$, $p \in \mathcal{P}_{2n}(\widetilde{A})$ e $p_n = \left( \begin{array}{ll} 1_n & 0 \\ 0 & 0 \end{array} \right)$ satisfazem
$$ \widetilde{\phi}(p) = v p_n v^*  \ \ \ \ \ \textrm{ e } \ \ \ \ \ \widetilde{\psi}(v) = \left( \begin{array}{ll} u & 0 \\ 0 & u^* \end{array} \right).$$
\end{Def}

Uma prova do resultado encutido na definição anterior pode ser encontrada em $9.1.4$ de \cite{RLL}. Além disso, o nome aplicação do índice se deve ao caso particular do índice de Fredholm, isto é,  quando a seqüência exata da definição anterior é dada por 
$$ 0 \rightarrow \mathcal{K} \hookrightarrow A \rightarrow A/\mathcal{K} \rightarrow 0$$
temos $\delta_1 : K_1(A/\mathcal{K}) \rightarrow K_0(\mathcal{K})$ como sendo o índice de Fredholm, vide capítulo $14$ de \cite{WO}.

\begin{Def} \label{SA}
Dada uma C$^*$-álgebra $A$ qualquer chamamos de suspensão de $A$ o conjunto simbolizado por $SA$ e dado por 
\begin{equation*}
\begin{split}
SA &= \{ f \in C([0,1],A): f(0) = f(1)= 0 \} \\
 &= \{ f \in C(S^1,A): f(1) = 0 \} = C_0(]0,1[, A),\\
\end{split}
\end{equation*}
e de cone de $A$ o conjunto denotado por $CA$ e dado por
$$ CA = \{ f \in C([0,1],A) : f(0) = 0\} .$$
\end{Def}

Através destas definições temos a seqüência exata curta
\begin{equation} \label{cone}
0 \rightarrow SA \overset{i}{\hookrightarrow} CA \overset{\pi}{\rightarrow} A \rightarrow 0
\end{equation}
para $\pi (f) = f(1)$.
E como $CA$ é homotopicamente equivalente a zero, temos $K_0(CA) = K_1(CA) = 0$ e assim a aplicação do índice definida anteriormente se torna um isomorfismo, como enunciamos a seguir.

\begin{Cor} \label{thetaa}
Dada uma C$^*$-álgebra $A$ qualquer a aplicação $\theta_A$ que é a aplicação $\delta_1$, dada  em \ref{indice}, da seqüências exata curta  $\ref{cone}$
$$\theta_A = \delta_1 : K_1(A) \rightarrow K_0(SA)$$ 
é um isomorfismo de grupos.
\end{Cor}

Um dos alicerces da K-teoria é a periodicidade de Bott e para entender este isomorfismo precisamos de alguns resultados e definições preliminares. 

Dada uma projeção $p \in \mathcal{P}_n(A)$, onde $A$ é uma C$^*$-álgebra com unidade, defi\-nimos $f_p \in C(S^1 , \mathcal{U}_n(A))$ dado por $f_p(z) = zp + (p_n - p)$.
Sabendo que neste caso a unitização de $SA$, isto é, $\widetilde{SA}= \{ (f, z) : f \in SA \textrm{ e } z \in \mathbb{C} \}$ pode ser identificado com o conjunto das aplicações $f \in C(S^1, A)$ tais que $f(1) \in \mathbb{C}1_A$, obtemos $f_p \in \mathcal{U}_n(\widetilde{SA})$. 

E assim definimos a aplicação de Bott para C$^*$-álgebras com unidade
\begin{equation*}
\begin{split}
\beta_A: K_0(A) &\rightarrow K_1(SA) \\
           [p]_0 &\rightarrow [f_p]_1.
\end{split}
\end{equation*}

Para C$^*$-álgebras sem unidade os resultados são análogos, porém como o retrato de $K_0(A)$ é um pouco diferente do anterior teremos a aplicação de Bott dada por 
\begin{equation*}
\begin{split}
\beta_A: K_0(A) &\rightarrow K_1(SA) \\
           [p]_0 - [s(p)]_0 &\rightarrow [f_pf^*_{s(p)}]_1.
\end{split}
\end{equation*}

Enunciaremos a seguir este teorema, cuja demonstração pode ser encontrada nas seções $9.1$ e $9.2$ de \cite{WO}, que prova que a aplicação de Bott definida acima é de fato um isomorfismo.

\begin{Teo}[Periodicidade de Bott] \label{Bott}
A aplicação $\beta_A: K_0(A) \rightarrow K_1(SA)$ é um isomorfismo de grupos.
\end{Teo}

Utilizando estes dois últimos resultados obtemos uma das principais ferramentas da K-teoria que é a seqüência exata cíclica de seis termos. Vejamos
 
\begin{Teo} \label{seis}
Dada a seqüência exata curta de C$^*$-álgebras
$$ 0 \rightarrow A \overset{\phi}{\rightarrow} B \overset{\psi}{\rightarrow} C \rightarrow 0$$
podemos construir a aplicação exponencial $\delta_0: K_0(C) \rightarrow K_1(A)$ que é dada pela  composição das aplicações
$$ K_0(C) \overset{\beta_C}{\rightarrow} K_1(SC) \simeq K_2(C) \overset{\delta_2}{\rightarrow} K_1(A) .$$ 
Disto e de \ref{indice} decorre a sequência exata cíclica 
$$ \xymatrix{
K_0(A)  \ar[rr]^-{K_0(\phi)}  & \ \ &  K_0(B) \ \ \ar[rr]^-{K_0(\psi)}     & \ \ &   K_0(C) \ar[dd]^-{\delta_0}     \\
& \\
K_1(C) \ar[uu]^-{\delta_1} & \ \ &  K_1(B)  \ar[ll]_-{K_1(\psi)}  & \ \  & K_1(A) \ar[ll]_-{K_1(\phi)}                     }. $$
\end{Teo}

De posse destas ferramentas e resultados podemos provar muitas outras propriedades importantes, como também exemplos bastante utilizados na literatura da área e neste trabalho.

\begin{Ex} \label{Ks1}
Para a C$^*$-álgebra $C(S^1)$, temos
$$K_0(C(S^1)) \simeq [\mathfrak{1}]_0 \mathbb{Z} \ \ \ \ \textrm{ e } \ \ \ \ K_1(C(S^1)) \simeq [\mathfrak{z}]_1\mathbb{Z}.$$
Onde $\mathfrak{1}$ simboliza a projeção $\mathfrak{1}(e^{i \theta}) = 1$ e $\mathfrak{z}$ o unitário $\mathfrak{z}(e^{i \theta}) = e^{i \theta}$.
\end{Ex}

Dado um espaço compacto Hausdorff $X$ e portanto a C$^*$-álgebra $C(S^1 \times X)$ sabemos que $C(S^1 \times X) \simeq C(S^1) \otimes X \simeq C(S^1, C(X)) $ e assim podemos construir a seqüência exata
$$ 0 \rightarrow SC(X) \rightarrow C(S^1, C(X)) \rightarrow C(X) \rightarrow 0$$
que cinde. Como $K_0$ e $K_1$ são funtores que preservam seqüências exatas cindidas, temos
$$ 0 \rightarrow K_i(SC(X)) \rightarrow K_i(C(S^1, C(X))) \rightarrow K_i(C(X)) \rightarrow 0$$
e portanto
\begin{Prop} \label{cs1a}
Nas condições anteriores temos para $i=0$ ou $1$ 
$$K_i(C(S^1, X)) \simeq K_i(X) \oplus K_{1-i}(X).$$
\end{Prop}

Além deste, alguns outros resultados sobre K-teoria são utilizados neste trabalho, mas estes, por aparecerem num contexto diferente do apresentado até agora, são enunciados e detalhados apenas quando utilizados o que pode não ser o mais didático, mas talvez seja o mais simples e prático.


\ \ \ 

\newpage

\ \ \

\end{document}